\documentclass{article}

\usepackage[utf8]{inputenc}
\usepackage[T1]{fontenc}
\usepackage{adjustbox}
\usepackage{amsmath}
\usepackage{amsthm}
\usepackage{array}
\usepackage{relsize}
\usepackage{mathpartir}
\usepackage{mathtools}
\usepackage{caption}
\usepackage{rotating}
\usepackage{tcolorbox}
\tcbuselibrary{
  skins,
  breakable
}
\usepackage{quiver}
\usepackage{titling} 
\usepackage{titlesec}
\usepackage{float}
\usepackage{graphics, graphicx, tikz}
\usepackage{tikz-cd}
\usetikzlibrary{
  arrows,
  arrows.meta,
  decorations.markings,
  shapes.geometric,
  trees
}
\usepackage{xspace}
\usepackage{xcolor}
\usepackage{geometry}
\usepackage[colorlinks=true, linktocpage=true]{hyperref}

\usepackage{biblatex}
\usepackage{mathpazo}
\usepackage[scaled=0.95]{helvet}
\usepackage{courier}
\titlespacing*{\subsection} {0pt}{0ex plus 1ex minus .2ex}{2ex plus .2ex}

\definecolor{grey}{rgb}{0.5,0.5,0.5}
\definecolor{mypurple}{rgb}{0.7,0,0.7}
\definecolor{links_color}{rgb}{0.4,0.1,0.1}
\hypersetup{
  linkcolor = links_color,
  citecolor = links_color
}

\title{A poset-like approach to positive opetopes}
\author{Louise \textsc{Leclerc}}

\counterwithin{figure}{section}

\newcounter{paragraphCounter}[section]
\renewcommand{\theparagraphCounter}{\thesection.\arabic{paragraphCounter}}

\newtheoremstyle{break}
  {}{}
  {}{}
  {\bfseries}{}
  {\newline}{}
\theoremstyle{break}

\newtheorem{rmk}[paragraphCounter]{}
\newtheorem{exm}[paragraphCounter]{}

\newenvironment{thm}[1][]{
  \refstepcounter{paragraphCounter}
  \begin{tcolorbox}
    [ enhanced
    , breakable
    , title =
        \textcolor{black}{
          \ifstrempty{#1}
            {Theorem~\theparagraphCounter}
            {Theorem~\theparagraphCounter\,:~#1}
        }
    , attach boxed title to top left = 
        { xshift = 15pt
        , yshift* = -\tcboxedtitleheight/2
        }
    , boxed title style = 
        { size = small
        , colback = purple!20
        , top = 1pt
        , bottom = 1pt
        }
    , colback = white
    , after skip = 10pt
    , before skip = 6pt
    ]
  }{\end{tcolorbox}}

\newenvironment{thm*}[1][]{
  \begin{tcolorbox}
    [ enhanced
    , breakable
    , title =
        \textcolor{black}{
          \ifstrempty{#1}
            {Theorem}
            {Theorem\,:~#1}
        }
    , attach boxed title to top left = 
        { xshift = 15pt
        , yshift* = -\tcboxedtitleheight/2
        }
    , boxed title style = 
        { size = small
        , colback = purple!20
        , top = 1pt
        , bottom = 1pt
        }
    , colback = white
    , after skip = 10pt
    , before skip = 6pt
    ]
  }{\end{tcolorbox}}

\newenvironment{defin}[1][]{
  \refstepcounter{paragraphCounter}
  \begin{tcolorbox}
    [ enhanced
    , breakable
    , title =
        \textcolor{black}{
          \ifstrempty{#1}
            {Definition~\theparagraphCounter}
            {Definition~\theparagraphCounter\,:~#1}
        }
    , attach boxed title to top left = 
        { xshift = 15pt
        , yshift* = -\tcboxedtitleheight/2
        }
    , boxed title style = 
        { size = small
        , colback = blue!30
        , top = 1pt
        , bottom = 1pt
        }
    , colback = white
    , after skip = 10pt
    , before skip = 6pt
    ]
  }{\end{tcolorbox}}

\newenvironment{prop}[1][]{
  \refstepcounter{paragraphCounter}
  \begin{tcolorbox}
    [ enhanced
    , breakable
    , title =
        \textcolor{black}{
          \ifstrempty{#1}
            {Proposition~\theparagraphCounter}
            {Proposition~\theparagraphCounter\,:~#1}
        }
    , attach boxed title to top left = 
        { xshift = 15pt
        , yshift* = -\tcboxedtitleheight/2
        }
    , boxed title style = 
        { size = small
        , colback = green!50!blue!30!
        , top = 1pt
        , bottom = 1pt
        }
    , colback = white
    , after skip = 10pt
    , before skip = 6pt
    ]
  }{\end{tcolorbox}}

\newenvironment{lem}[1][]{
  \refstepcounter{paragraphCounter}
  \begin{tcolorbox}
    [ enhanced
    , breakable
    , title =
        \textcolor{black}{
          \ifstrempty{#1}
            {Lemma~\theparagraphCounter}
            {Lemma~\theparagraphCounter\,:~#1}
        }
    , attach boxed title to top left = 
        { xshift = 15pt
        , yshift* = -\tcboxedtitleheight/2
        }
    , boxed title style = 
        { size = small
        , colback = blue!30
        , top = 1pt
        , bottom = 1pt
        }
    , colback = white
    , after skip = 10pt
    , before skip = 6pt
    ]
  }{\end{tcolorbox}}

\newcommand{\para}{
  \refstepcounter{paragraphCounter}
  \noindent\bfseries{\theparagraphCounter.}
  \normalfont
}

\newcommand{\refrmk}[1]{\hyperref[#1]{Remark \ref*{#1}}}
\newcommand{\refexm}[1]{\hyperref[#1]{Exemple \ref*{#1}}}
\newcommand{\refthm}[1]{\hyperref[#1]{Theorem \ref*{#1}}}
\newcommand{\refdefin}[1]{\hyperref[#1]{Definition \ref*{#1}}}
\newcommand{\refprop}[1]{\hyperref[#1]{Proposition \ref*{#1}}}
\newcommand{\reflem}[1]{\hyperref[#1]{Lemma \ref*{#1}}}

\newcommand{\tailleftline}{\DOTSB\relbar\joinrel\hspace{-4px}\raisebox{0.3px}{\mbox{\small<}}}

\newcommand{\tailrightline}{\,\raisebox{0.3px}{\mbox{\small>}}\hspace{-2.5px}\joinrel\DOTSB\relbar}
\newcommand{\branch}[1]{\ensuremath{\tailleftline_{#1}\,}}

\newcommand{\branched}[1]{\ensuremath{\tailrightline_{#1}}}

\newcommand{\lbl}[1]{\ensuremath{\scriptstyle #1}}
\newcommand{\N}{\ensuremath{\mathbb{N}}\xspace}
\newcommand{\cl}[1]{\ensuremath{\mathsf{cl}\!\left(#1\right)}}
\newcommand*{\intInter}[2]
  {\ensuremath{
    \left[\!\left[\, \ensuremath{{#1}},\, \ensuremath{{#2}} \,\right]\!\right]
  }}
\renewcommand{\dim}{\ensuremath{\mathsf{dim}}}

\newcommand{\trigup}[1]{\ensuremath{\triangleleft^{\scriptscriptstyle \,C_{#1},\,+}}}

\begin{document}

  \maketitle

  
\begin{abstract}
  We introduce in this paper a new formalisation of positive opetopes where faces are organised in a poset.
  Then we show that our definition is equivalent to that of positives opetopes as given by Marek \textsc{Zawadowski} in \cite{zawadowski2023positive}.
\end{abstract}

  
\section*{Introduction}
\label{sec:intro}

We begin by motivating the notion of opetope that is lying at the heart of this article.
When manipulating categories, we find ourselves considering two types of forms:
zero dimensional elements, which are the objects (the $0$-cells):
$$\bullet$$
And one dimensional ones: the arrows (the $1$-cells).
$$
\bullet \longrightarrow \bullet
$$
Category theory is mainly about composing arrows, so the category-theorist often draws diagrams as below:
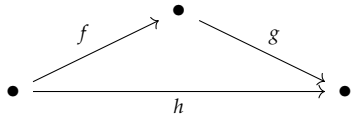
\[\begin{tikzcd}
  && \bullet \\
  \bullet &&&& \bullet
  \arrow["f", from=2-1, to=1-3]
  \arrow["g", from=1-3, to=2-5]
  \arrow["h"', from=2-1, to=2-5]
\end{tikzcd}\]
and says that it commutes if $g \circ f = h$.
It can be understood as a kind of surface whose boundary is given by the arrows $f,\,g$ and $h$, which ensures that there is a way to go from the upper path to the lower one.
Since in a weak setting we do not like equalities, it would be better to replace this $g \circ f = h$ by some kind of algebraic data $$\textcolor{blue}{\alpha} : g \circ f \,\textcolor{blue}{\Longrightarrow}\, h$$
Hence we give a name to our (oriented) surface (or $2$-cell), $\textcolor{blue}{\alpha}$, and depict it as below:
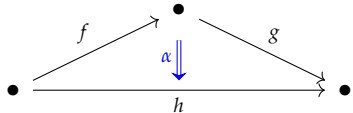
\[\begin{tikzcd}[
  two/.style = {-{Stealth[length=3pt, width=6pt]}, double, double distance = 1pt, draw = blue, shorten >= 5pt, shorten <= 5pt}
  ]
  && \bullet \\
  \bullet &&&& \bullet
  \arrow["f", from=2-1, to=1-3]
  \arrow["g", from=1-3, to=2-5]
  \arrow[""{name=0, anchor=center, inner sep=0}, "h"', from=2-1, to=2-5]
  \arrow["\alpha"'{blue}, two, shorten >=3pt, from=1-3, to=0]
\end{tikzcd}\]
More generally, one is interested in $n$-ary compositions of arrows, so we may also consider diagrams such as
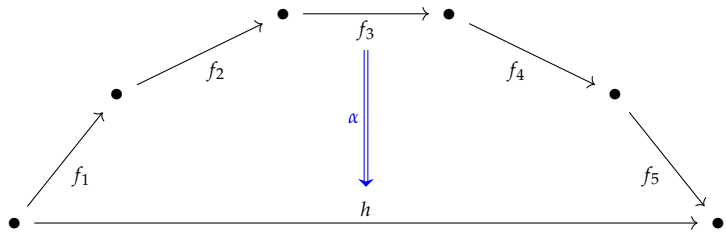
\[\begin{tikzcd}[
    two/.style = {-{Stealth[length=3pt, width=6pt]}, double, double distance = 1pt, draw = blue, shorten >= 5pt, shorten <= 5pt}
  ]
  &&& \bullet && \bullet \\
  & \bullet &&&&&& \bullet \\
  \\
  \bullet &&&&&&&& \bullet
  \arrow["{f_1}"', from=4-1, to=2-2]
  \arrow["{f_2}"', from=2-2, to=1-4]
  \arrow[""{name=0, anchor=center, inner sep=0}, "{f_3}"', from=1-4, to=1-6]
  \arrow["{f_4}"', from=1-6, to=2-8]
  \arrow["{f_5}"', from=2-8, to=4-9]
  \arrow[""{name=1, anchor=center, inner sep=0}, "h", from=4-1, to=4-9]
  \arrow["\alpha"'{blue}, two, shorten <=13pt, shorten >=13pt, from=0, to=1]
\end{tikzcd}\]
Since our surfaces are now algebraic data, we will soon be interested in comparing them, too.
Once again, since asking if they are equal or not is more a strict-minded approach to category theory, we wish to introduce named volumes ($3$-cells), which may assert that there is a way to go from a composite of surfaces to another one. An illustration is given below:
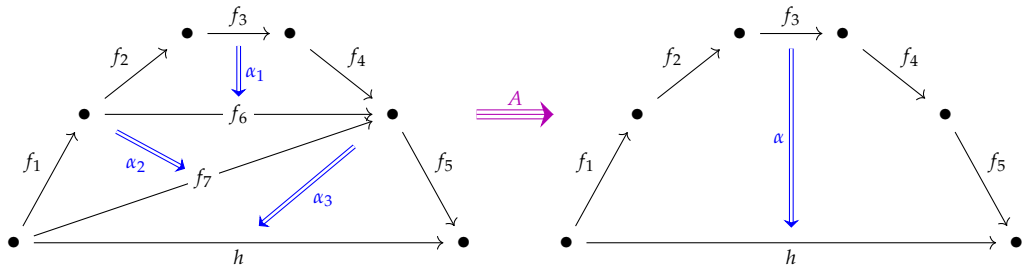
\[\begin{tikzcd}[
    column sep = small,
  , two/.style = {-{Stealth[length=3pt, width=6pt]}, double, double distance = 1pt, draw = blue, shorten >= 5pt, shorten <= 5pt}
  , three/.style = {preaction={draw = mypurple, -{Stealth[length=5pt, width=10pt]}, double, double distance = 3pt}, draw = mypurple, -}
  ]
  &&& \bullet && \bullet &&&&&&&& \bullet && \bullet \\
  & \bullet &&&&&& \bullet & {} && {} & \bullet &&&&&& \bullet \\
  \\
  \bullet &&&&&&&& \bullet && \bullet &&&&&&&& \bullet
  \arrow["{f_1}", from=4-1, to=2-2]
  \arrow["{f_2}", from=2-2, to=1-4]
  \arrow[""{name=0, anchor=center, inner sep=0}, "{f_3}", from=1-4, to=1-6]
  \arrow["{f_4}", from=1-6, to=2-8]
  \arrow["{f_5}", from=2-8, to=4-9]
  \arrow[""{name=1, anchor=center, inner sep=0}, "h"', from=4-1, to=4-9]
  \arrow[""{name=2, anchor=center, inner sep=0}, "{f_6}"{description}, from=2-2, to=2-8]
  \arrow[""{name=3, anchor=center, inner sep=0}, "{f_7}"{description}, from=4-1, to=2-8]
  \arrow["A"{mypurple, above = 1pt}, three, from=2-9, to=2-11]
  \arrow["{f_1}", from=4-11, to=2-12]
  \arrow["{f_2}", from=2-12, to=1-14]
  \arrow[""{name=4, anchor=center, inner sep=0}, "{f_3}", from=1-14, to=1-16]
  \arrow["{f_4}", from=1-16, to=2-18]
  \arrow["{f_5}", from=2-18, to=4-19]
  \arrow[""{name=5, anchor=center, inner sep=0}, "h"', from=4-11, to=4-19]
  \arrow["{\alpha_1}"{blue}, two, shorten <=4pt, shorten >=6pt, from=0, to=2]
  \arrow["{\alpha_2}"'{blue}, two, shorten <=5pt, shorten >=7pt, from=2-2, to=3]
  \arrow["{\alpha_3}"{blue}, two, shorten <=9pt, shorten >=9pt, from=2-8, to=1]
  \arrow["\alpha"'{blue}, two, shorten <=5pt, shorten >=5pt, from=4, to=5]
\end{tikzcd}\]
And as another example:
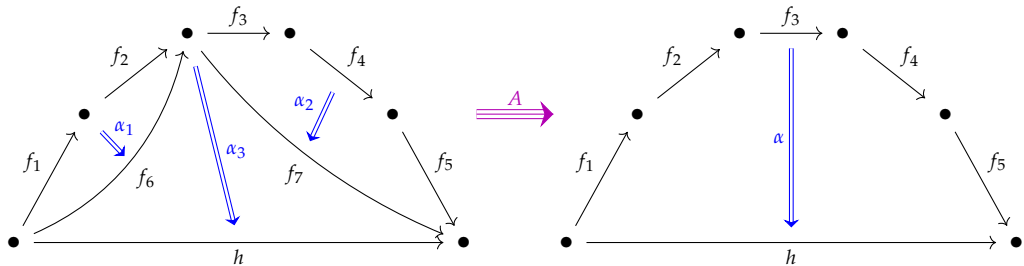
\[\begin{tikzcd}[
    column sep = small,
  , two/.style = {-{Stealth[length=3pt, width=6pt]}, double, double distance = 1pt, draw = blue, shorten >= 5pt, shorten <= 5pt}
  , three/.style = {preaction={draw = mypurple, -{Stealth[length=5pt, width=10pt]}, double, double distance = 3pt}, draw = mypurple, -}
  ]
  &&& \bullet && \bullet &&&&&&&& \bullet && \bullet \\
  & \bullet &&&&&& \bullet & {} && {} & \bullet &&&&&& \bullet \\
  \\
  \bullet &&&&&&&& \bullet && \bullet &&&&&&&& \bullet
  \arrow["{f_1}", from=4-1, to=2-2]
  \arrow["{f_2}", from=2-2, to=1-4]
  \arrow["{f_3}", from=1-4, to=1-6]
  \arrow[""{name=0, anchor=center, inner sep=0}, "{f_4}", from=1-6, to=2-8]
  \arrow["{f_5}", from=2-8, to=4-9]
  \arrow[""{name=1, anchor=center, inner sep=0}, "h"', from=4-1, to=4-9]
  \arrow["A"{mypurple, above = 1pt}, three, from=2-9, to=2-11]
  \arrow["{f_1}", from=4-11, to=2-12]
  \arrow["{f_2}", from=2-12, to=1-14]
  \arrow[""{name=2, anchor=center, inner sep=0}, "{f_3}", from=1-14, to=1-16]
  \arrow["{f_4}", from=1-16, to=2-18]
  \arrow["{f_5}", from=2-18, to=4-19]
  \arrow[""{name=3, anchor=center, inner sep=0}, "h"', from=4-11, to=4-19]
  \arrow["{f_6}"'{name=4}, curve={height=18pt}, from=4-1, to=1-4]
  \arrow["{f_7}"'{name=5}, curve={height=12pt}, from=1-4, to=4-9]
  \arrow["\alpha"'{blue}, two, shorten <=5pt, shorten >=5pt, from=2, to=3]
  \arrow["{\alpha_1}"{blue, above = 7pt, right = -3pt}, two, shorten >=3pt, shorten <=0pt, from=2-2, to=4]
  \arrow["{\alpha_2}"'{blue}, two, shorten <=7pt, shorten >=7pt, from=0, to=5]
  \arrow["{\alpha_3}"{blue, pos = 0.6}, two, shorten <=6pt, shorten >=6pt, from=1-4, to=1]
\end{tikzcd}\]
Typically, in those drawings, (which are in fact $3$-dimensional \emph{opetopes}) the rightmost $2$-cell $\textcolor{blue}{\alpha}$ will be called the \emph{target} of the $3$-cell $\textcolor{mypurple}{A}$.
The $2$-cells $\textcolor{blue}{\alpha_1}$, $\textcolor{blue}{\alpha_2}$ and $\textcolor{blue}{\alpha_3}$ will be called its \emph{sources}.
Notice that each $1$-cell (our basic arrows $\bullet \to \bullet$) has exacly one source, and one target.
But in higher dimensions, $2$-cells and $3$-cells may have many sources, while still having exactly one target.
We may also say that $0$-dimensional cells, (\textit{i.e.} points $\bullet$), have no sources, and no target.
Pictures as above take the shape of a polytope when increasing dimension, and they represent algebraic operations, whence the term introduced by \textsc{Baez} and \textsc{Dolan} in \cite{baez1997higherdimensional}: \emph{opetopes}, for "ope(ration-poly)topes".

We have several ways to encode such shapes combinatorially.
\begin{itemize}
  \item An approach, taken by Marek \textsc{Zawadowski} in \cite{zawadowski2023positive}, is to name all cells as above (including 0-cells), and store them in a poset, where every element has a dimension (its geometric dimension), and a relation $z \le x$ will mean that $z$ is a subface of $x$ in the opetope.
    Then we need to identify axioms to ensure that those posets fit our intuition of opetopes as above. This is exposed in Section \ref{sec:ZPO}.

    Using the ideas of Amar \textsc{Hadzihasanovic} \cite{hadzihasanovic2019combinatorialtopological} and Pierre-Louis \textsc{Curien}, we were able to identify a second formalism using a similar principle, presented in Section \ref{sec:DFC}.

    In the sections \ref{sec:dfc_to_zpo} and \ref{sec:zpo_to_dfc}, we will show that the two formalisms are equivalent.
  \item We may also notice that there are trees hidden in opetopes, which may be retrieved by \textsc{Poincaré} duality:
    
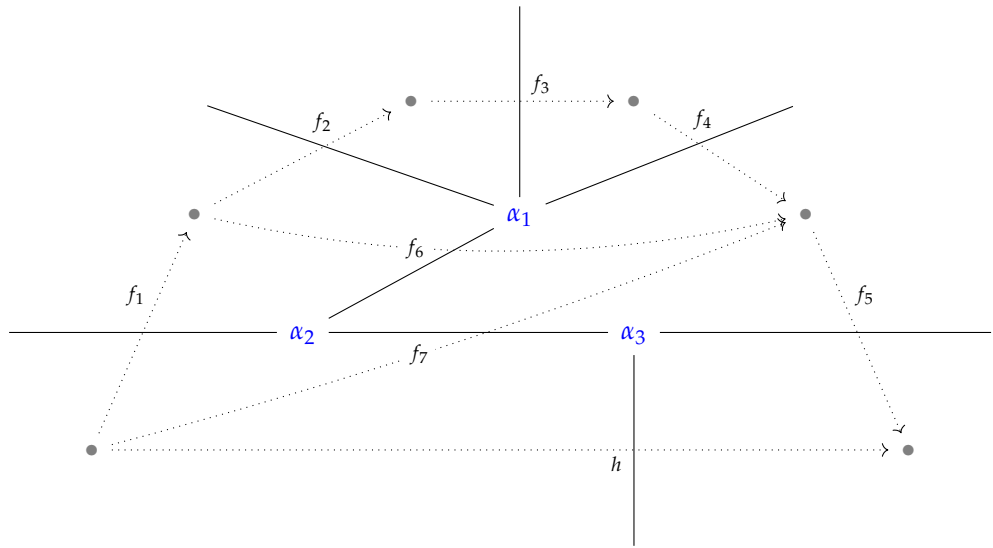
\vspace{-5pt}
    \[\begin{tikzcd}[row sep = 30pt]
      &&&&& {} \\
      && {} && {\textcolor{grey}{\bullet}} && {\textcolor{grey}{\bullet}} && {} \\
      && {\textcolor{grey}{\bullet}} &&& {\textcolor{blue}{\alpha_1}} &&& {\textcolor{grey}{\bullet}} \\
      {} &&& {\textcolor{blue}{\alpha_2}} &&& {\textcolor{blue}{\alpha_3}} &&&& {} \\
      & {\textcolor{grey}{\bullet}} &&&&&&&& {\textcolor{grey}{\bullet}} \\
      &&&&&& {}
      \arrow["{f_1}"{pos=0.6}, dotted, from=5-2, to=3-3]
      \arrow["{f_2}"{pos=0.7}, dotted, from=3-3, to=2-5]
      \arrow["{f_3}"{pos=0.6}, dotted, from=2-5, to=2-7]
      \arrow["{f_4}"{pos=0.25}, dotted, from=2-7, to=3-9]
      \arrow["{f_5}"{pos=0.4}, dotted, from=3-9, to=5-10]
      \arrow["h"'{pos=0.65}, dotted, from=5-2, to=5-10]
      \arrow["{f_6}"{description, pos=0.35}, curve={height=18pt}, dotted, from=3-3, to=3-9]
      \arrow["{f_7}"{description, pos=0.45}, curve={height=6pt}, dotted, from=5-2, to=3-9]
      \arrow[no head, from=6-7, to=4-7]
      \arrow[no head, from=4-7, to=4-4]
      \arrow[no head, from=4-4, to=3-6]
      \arrow[no head, from=3-6, to=1-6]
      \arrow[no head, from=3-6, to=2-3]
      \arrow[no head, from=3-6, to=2-9]
      \arrow[no head, from=4-4, to=4-1]
      \arrow[no head, from=4-7, to=4-11]
    \end{tikzcd}\]
    And for the second example:
    
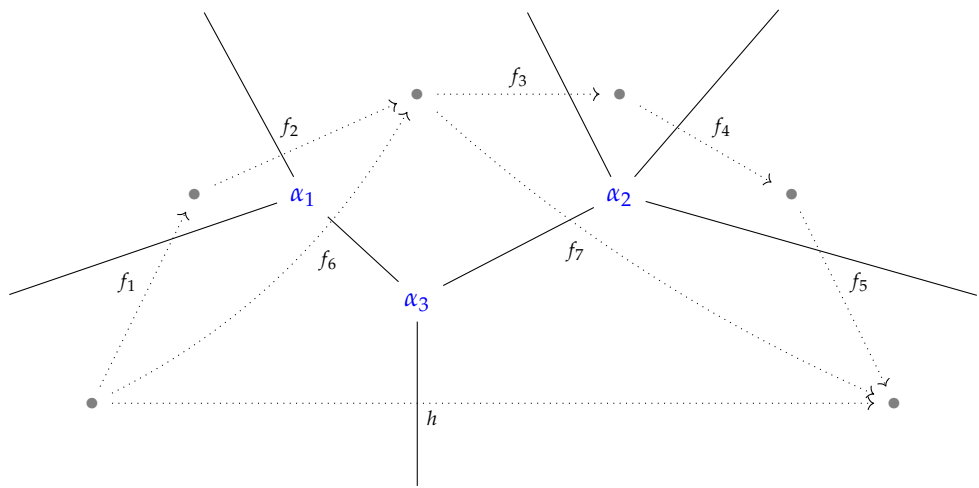
\[\begin{tikzcd}[row sep = 25pt]
      && {} &&& {} &&& {} \\
      &&&& {\textcolor{grey}{\bullet}} && {\textcolor{grey}{\bullet}} \\
      && {\textcolor{grey}{\bullet}} & {\textcolor{blue}{\alpha_1}} &&& {\textcolor{blue}{\alpha_2}} && {\textcolor{grey}{\bullet}} \\
      {} &&&& {\textcolor{blue}{\alpha_3}} &&&&&& {} \\
      & {\textcolor{grey}{\bullet}} &&&&&&&& {\textcolor{grey}{\bullet}} \\
      &&&& {}
      \arrow["{f_1}", dotted, from=5-2, to=3-3]
      \arrow["{f_2}", dotted, from=3-3, to=2-5]
      \arrow["{f_3}", dotted, from=2-5, to=2-7]
      \arrow["{f_4}", dotted, from=2-7, to=3-9]
      \arrow["{f_5}", dotted, from=3-9, to=5-10]
      \arrow["h"'{pos=0.42}, dotted, from=5-2, to=5-10]
      \arrow["{f_6}"'{pos=0.61}, curve={height=18pt}, dotted, from=5-2, to=2-5]
      \arrow["{f_7}"'{pos=0.38}, curve={height=12pt}, dotted, from=2-5, to=5-10]
      \arrow[no head, from=6-5, to=4-5]
      \arrow[no head, from=4-5, to=3-4]
      \arrow[no head, from=3-4, to=1-3]
      \arrow[no head, from=3-4, to=4-1]
      \arrow[no head, from=4-5, to=3-7]
      \arrow[no head, from=3-7, to=1-6]
      \arrow[no head, from=3-7, to=1-9]
      \arrow[no head, from=3-7, to=4-11]
    \end{tikzcd}\]
\end{itemize}
Here, cells $\textcolor{blue}{\alpha_1}$, $\textcolor{blue}{\alpha_2}$ and $\textcolor{blue}{\alpha_3}$ are lower dimensional opetopes, and should have an associated tree too.
Hence we may represent an opetope with a bunch of trees in several dimensions, interconnected by gluing relations.
This is the approach taken by Joachim \textsc{Kock}, André \textsc{Joyal}, Michael \textsc{Batanin} and Jean-François \textsc{Mascari} in \cite{Kock2010}, or by Cédric \textsc{Ho Thanh}, Pierre-Louis \textsc{Curien} and Samuel \textsc{Mimram} in \cite{Curien2022} 
And also, from a different perspective originating from Pierre-Louis \textsc{Curien}, with the formalisation of \emph{epiphytes}, which will be presented in a future paper.

In order to provide the reader with a better intuition of opetopes in higher dimension, an illustration of a $4$-dimensional opetope is given below:
\vspace{15pt}



%
  \adjustbox{scale=0.9, center}
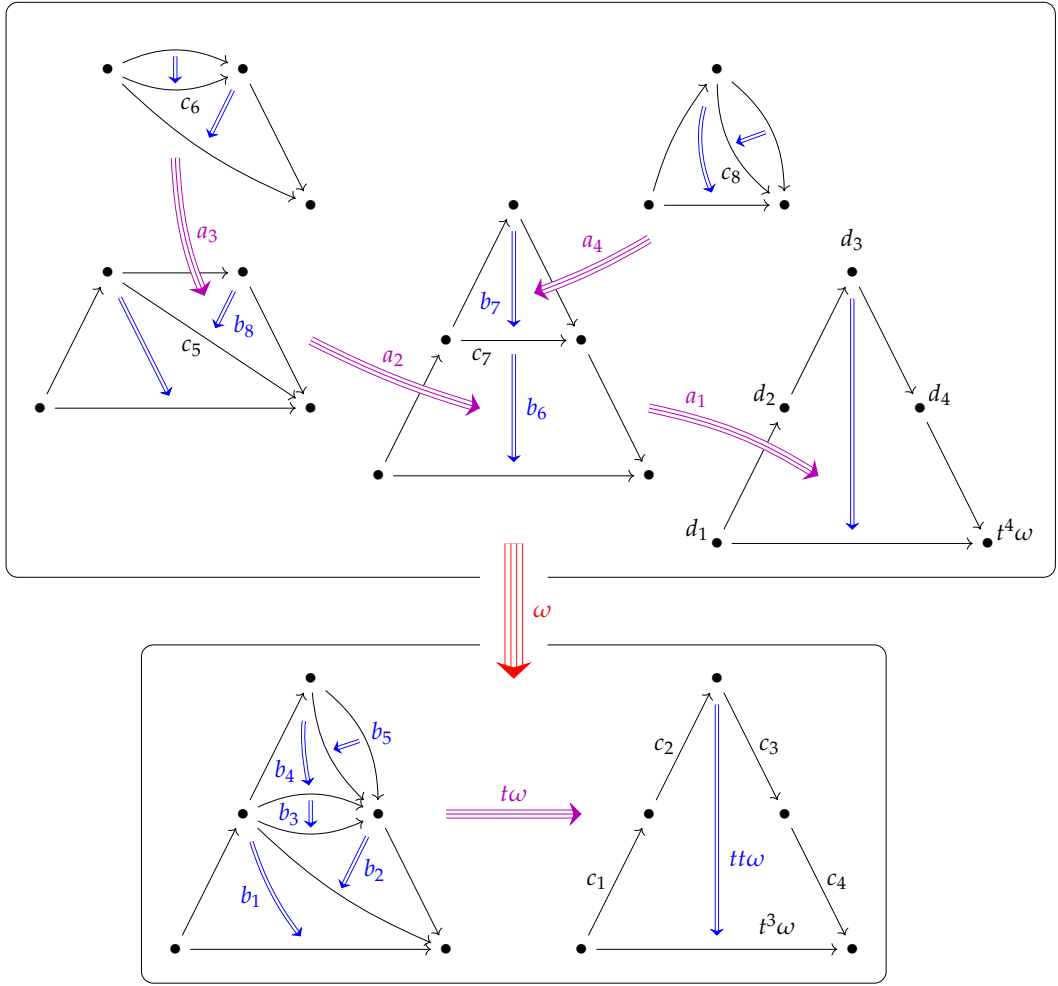
{
  \begin{tikzpicture}
    [
    , along/.style = {midway, inner sep = 0pt}
    , two/.style = {-{Stealth[length=3pt, width=6pt]}, double, double distance = 1pt, draw = blue, shorten >= 5pt, shorten <= 5pt}
    , three/.style={preaction={draw = mypurple, -{Stealth[length=5pt, width=10pt]}, double, double distance = 3pt}, draw = mypurple, -}
    , four/.style={preaction={draw = red, -{Stealth[length=7pt, width=15pt]}, double, double distance = 7pt, shorten >=0pt}, draw = red, -, double, double distance=2pt, shorten >=5pt}
    ]

    \node (d1_a1)  at (-2, 0)  {$\bullet$};
    \node (d2_a1)  at (-1, 2)  {$\bullet$}
      edge [<-] (d1_a1);
    \node (d3_a1)  at (0, 4)   {$\bullet$}
      edge [<-] (d2_a1);
    \node (d4_a1)  at (1, 2)   {$\bullet$}
      edge [<-] (d3_a1)
      edge [<-]
        node[below, pos = 0.8] {$c_7$}
        node[along] (c7_a1) {}
        (d2_a1);
    \node (d5_a1)  at (2, 0)   {$\bullet$}
      edge [<-] (d4_a1)
      edge [<-]
        node[along] (ttt_a1) {} 
        (d1_a1);
    \draw[two] (d3_a1) to
      node[pos = 0.7, left, color = blue, inner sep = 5pt] {$b_7$}
      (c7_a1);
    \draw[two] (c7_a1) to
      node[right, color = blue, inner sep = 5pt] {$b_6$} (ttt_a1);

    \node (d1_a2) at (-7, 1)  {$\bullet$};
    \node (d2_a2) at (-6, 3)  {$\bullet$}
      edge [<-] (d1_a2);
    \node (d4_a2) at (-4, 3)  {$\bullet$}
      edge [<-] (d2_a2);
    \node (d5_a2) at (-3, 1)  {$\bullet$}
      edge [<-] (d4_a2)
      edge [<-]
        node[below = 1pt, pos = 0.6] {$c_5$}
        node[along] (c5_a2) {}
        (d2_a2)
      edge [<-]
        node[along] (ttt_a2) {}
        (d1_a2);
    \draw[two, shorten <= 2pt] (d4_a2) to
      node[pos = 0.75, right, color = blue, inner sep = 7pt] {$b_8$} (c5_a2);
    \draw[two] (d2_a2) -- (ttt_a2);

    \node (d2_a3) at (-6, 6)  {$\bullet$};
    \node (d4_a3) at (-4, 6)  {$\bullet$}
      edge [<-, bend right = 25]
        node[along] (c7_a3) {}
        (d2_a3)
      edge [<-, bend left = 25] 
        node[below, pos = 0.35] {$c_6$}
        node[along] (c6_a3) {}
        (d2_a3);
    \node (d5_a3) at (-3, 4)  {$\bullet$}
      edge [<-] (d4_a3)
      edge [<-, bend left = 10] 
        node[along] (c5_a3) {}
        (d2_a3);
    \draw[two, shorten <= 2pt, shorten >= 2pt] (c7_a3) -- (c6_a3);
    \draw[two, shorten <= 3pt] (d4_a3) -- (c5_a3);

    \node (d2_a4) at (2, 4)   {$\bullet$};
    \node (d3_a4) at (3, 6)   {$\bullet$}
      edge [<-, bend right = 10] (d2_a4);
    \node (d4_a4) at (4, 4)   {$\bullet$}
      edge [<-, bend right = 25] 
        node[along] (c3_a4) {}
        (d3_a4)
      edge [<-, bend left = 25]
        node[left, pos = 0.2] {$c_8$}  
        node[along] (c8_a4) {}
        (d3_a4)
      edge [<-] 
        node[along] (c7_a4) {}
        (d2_a4);
    \draw[two, shorten <= 2pt, shorten >= 2pt] (c3_a4) -- (c8_a4);
    \draw[two, shorten <= 10pt] (d3_a4) to[bend right = 20] (c7_a4);

    \node[label = {[left]$d_1$}] (d1_ta1)  at (3, -1)  {$\bullet$};
    \node[label = {[left]$d_2$}] (d2_ta1)  at (4, 1)  {$\bullet$}
      edge [<-] (d1_ta1);
    \node[label = {[above]$d_3$}] (d3_ta1)  at (5, 3)   {$\bullet$}
      edge [<-] (d2_ta1);
    \node[label = {[right]$d_4$}] (d4_ta1)  at (6, 1)   {$\bullet$}
      edge [<-] (d3_ta1);
    \node[label = {[right]$t^4\omega$}] (d5_ta1)  at (7, -1)   {$\bullet$}
      edge [<-] (d4_ta1)
      edge [<-]
        node[along] (ttt_ta1) {}
        (d1_ta1);
    \draw[two] (d3_ta1) to (ttt_ta1);

    \node (d1_t)  at (-5, -7)  {$\bullet$};
    \node (d2_t)  at (-4, -5)  {$\bullet$}
      edge [<-] (d1_t);
    \node (d3_t)  at (-3, -3)   {$\bullet$}
      edge [<-] (d2_t);
    \node (d4_t)  at (-2, -5)   {$\bullet$}
      edge [<-, bend right = 25]
        node[along] (c9_t) {} 
        (d3_t)
      edge [<-, bend left = 20]
        node[along] (c8_t) {}
        (d3_t)
      edge [<-, bend right = 25]
        node[along] (c7_t) {}
        (d2_t)
      edge [<-, bend left = 25]
        node[along] (c6_t) {}
        (d2_t);
    \node (d5_t)  at (-1, -7)   {$\bullet$}
      edge [<-] (d4_t)
      edge [<-, bend left = 10]
        node[along] (c5_t) {}
        (d2_t)
      edge [<-]
        node[along] (ttt_t) {} 
        (d1_t);
    \draw[two, shorten <= 2pt, shorten >= 2pt] (c9_t) to
      node[right = 15pt, above = -3pt, color = blue] {$b_5$}
      (c8_t);
    \draw[two, shorten <= 12pt, shorten >= 3pt] (d3_t) to[bend right = 10] 
      node[pos = 0.8, left = 1pt, color = blue] {$b_4$}
      (c7_t);
    \draw[two, shorten <= 2pt, shorten >= 2pt] (c7_t) to 
      node[left = 1pt, color = blue] {$b_3$}
      (c6_t);
    \draw[two, shorten <= 4pt, shorten >= 2pt] (d4_t) to
      node[pos = 0.65, right = 2pt, color = blue] {$b_2$}
      (c5_t);
    \draw[two, shorten >= 6pt, shorten <= 6pt] (d2_t) to[bend right = 10] 
      node[pos = 0.55, left = 2pt, color = blue] {$b_1$}
      (ttt_t);

    \node (d1_tt)  at (1, -7)  {$\bullet$};
    \node (d2_tt)  at (2, -5)  {$\bullet$}
      edge [<-]
        node[left] {$c_1$}
        (d1_tt);
    \node (d3_tt)  at (3, -3)   {$\bullet$}
      edge [<-] 
        node[left] {$c_2$}
        (d2_tt);
    \node (d4_tt)  at (4, -5)   {$\bullet$}
      edge [<-]
        node[right] {$c_3$}
        (d3_tt);
    \node (d5_tt)  at (5, -7)   {$\bullet$}
      edge [<-]
        node[right] {$c_4$}
        (d4_tt)
      edge [<-]
        node[pos = 0.25, above, inner sep = 6pt] {$t^3\omega$}
        node[along] (ttt_tt) {}
        (d1_tt);
    \draw[two] (d3_tt) to
      node[pos = 0.65,right, color = blue, inner sep = 5pt] {$tt\omega$} 
      (ttt_tt);


    \draw[three] (-3, 2) to[bend right = 5]
      node[above, color = mypurple, inner sep = 5pt] {$a_2$}
      (-0.5,1);
    \draw[three, shorten >= 5pt] (-5, 4.7) to[bend right = 10] 
      node[right, color = mypurple, inner sep = 5pt] {$a_3$}
      (-4.5, 2.5);
    \draw[three] (2, 3.5) to[bend left = 5] 
      node[above, color = mypurple, inner sep = 7pt] {$a_4$}
      (0.3,2.7);
    \draw[three] (2, 1) to[bend left = 10] 
      node[near start, above, color = mypurple, inner sep = 5pt] {$a_1$}
      (4.5, 0);
    \draw[three] (-1, -5) to
      node[above, color = mypurple, inner sep = 5pt] {$t\omega$}
      (1, -5);


    \pgfsetcornersarced{\pgfpoint{5pt}{5pt}}
    \draw (-0.5, -1.5) -- (-7.5, -1.5) -- (-7.5, 7) -- (8, 7) -- (8, -1.5) -- (0.5, -1.5);
    \draw (-0.5, -2.5) -- (-5.5, -2.5) -- (-5.5, -7.5) -- (5.5, -7.5) -- (5.5, -2.5) -- (0.5, -2.5);
    \pgfsetcornersarced{\pgfpointorigin}

    \draw[four] (0, -1) to
      node[right, color = red, inner sep = 8pt] {$\omega$}
      (0, -3);

  \end{tikzpicture}  
  }
%

For the reader introduced to the ideas of higher category theory, we specify that opetopes are a form of cell, like globes, or simplexes. 
They fit into the broader context of opetopic cardinals, which plays a role similar to that of pasting diagrams (see \cite{leinster2003higher}, or \cite{Ara2010} where they are called \emph{globular sums}) in the globular setting.
In particular, opetopic cardinals are arranged in a strict $\omega$-category, freely generated by opetopes of all dimensions.
The study of opetopes is (for instance) motivated by the following result (see Corollary 13.5 in \cite{zawadowski2023positive}, which is proved by using the aforementioned opetopic cardinals):

\begin{thm*}[\textsc{Zawadowski}]
  There is an equivalence of categories between $\widehat{\mathbf{pOpe}}$ and $\mathbf{pPoly}$.
\end{thm*}

where $\mathbf{pOpe}$ denotes the category of positive opetopes, $\widehat{\mathbf{pOpe}}$ the associated presheaf category, and $\mathbf{pPoly}$ the category of positive-to-one polygraphs.
That is, positive opetopes may be used to present a strict $\omega$-category.

  
\section{Dendritic face complexes}
\label{sec:DFC}

  The poset structure that will be introduced below makes implicit use of a notion of rooted tree, which we start by defining below.
  \begin{defin}[Rooted tree]\label{defin:rooted_tree}
    A \emph{rooted tree} $T$ consist of:
    \begin{itemize}
      \item A \emph{finite} set of nodes $T^\bullet$.
      \item For each node $a \in T^\bullet$, a \emph{finite} set $A(a)$, called the arity of $a$.
      \item A (necessarily finite) set of \emph{triplets}, denoted $a \branch{b} a'$ for some $a,\,a' \in T^\bullet$ and $b\in A(a)$.
        Moreover we ask that for each $a\in T^\bullet$ and $b\in A(a)$, there is at most one triplet $a \branch{b} a'$.
        If there is none, the pair $(a,\,b)$ is said to be a \emph{leaf} of $T$, and we let
        \[ 
        T^\vert := \{ (a,\,b) \mbox{ leaf of } T \}
        \]
    \end{itemize}
    We moreover ask for a distinguished element $\rho(T) \in T^\bullet$, called the \emph{root} of $T$, satisfying the following property:
    for each node $a\in T^\bullet$, there is a unique \emph{(descending) path} in $T$
    \[
    a = a_0 \branched{b_1} a_1 \branched{b_1} \cdots \branched{b_p} a_p = \rho(T)
    \]
    from $a$ to the root of $T$.
  \end{defin}
  \begin{rmk}\label{rmk:epi_unique_root}
    Notice that if it exists, the root is uniquely determined.
  \end{rmk}
  \begin{exm}\label{exm:rooted_tree}
  Below is a representation of the rooted tree $T$ having
  \begin{itemize}
    \item as nodes $T^\bullet := \{a_1,\,a_2,\,a_3,\,a_4\}$
    \item as arities
      \begin{mathpar}
        A(a_1) := \{b_6,\,b_7\}
        \and
        A(a_2) := \{b_1,\,b_8\}
        \and
        A(a_3) := \{b_2,\,b_3\}
        \and
        A(a_4) := \{b_4,\,b_5\}
      \end{mathpar}
    \item as triplets
      \begin{mathpar}
        a_1 \branch{b_6} a_2
        \and
        a_1 \branch{b_7} a_4
        \and
        a_2 \branch{b_8} a_3
      \end{mathpar}
    \item as root $a_1$.
  \end{itemize}
  
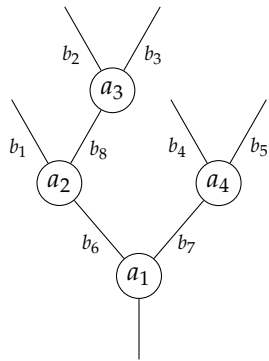
\[\begin{tikzpicture}
    [
      level distance = 35px
    , inner/.style = {circle, draw, inner sep=2px}
    , edge from parent/.style = {draw, near start, inner sep = 4px}
    , level 2/.style = {sibling distance=60px}
    , level 3/.style = {sibling distance=40px}
    , level 4/.style = {sibling distance=40px}
    ]
    \node[draw = none] {} [grow'=up] 
      child {
        node[inner] {$a_1$}
          child {
            node[inner] {$a_2$}
              child {node {} edge from parent node[left] {\lbl{b_1}}}
              child {
                node[inner] {$a_3$}
                  child {node {} edge from parent node[left] {\lbl{b_2}}}
                  child {node {} edge from parent node[right] {\lbl{b_3}}}
                edge from parent node[right] {\lbl{b_8}}
              }
            edge from parent node[left] {\lbl{b_6}}
          }
          child {
            node[inner] {$a_4$}
              child {node {} edge from parent node[left] {\lbl{b_4}}}
              child {node {} edge from parent node[right] {\lbl{b_5}}}
            edge from parent node[right] {\lbl{b_7}}
          }
      };
  \end{tikzpicture}\]
  \end{exm}

  \para We now present the notion of dendritic face complexes (DFC), defined in order to encode opetopes.
  Formally, they are partially ordered sets, where relations between elements are labelled by a sign ($+$ or $-$) in the Hasse diagram of which.
  They shall also satisfy some properties that will be given below.
  The elements of the poset, sometimes called \emph{faces}, stand for the faces of the opetope.
  All faces have a dimension in the poset, which reflects the intuitive dimension in the geometrical sense.
  A relation $y \prec x$ should be thought as "$y$ is a codimension-1 subface of $x$".
  Since in an opetope, faces come with an orientation, the sign tells if $y$ is an input or an output of $x$.

  For exemple, if the picture belows is part of an opetope, we should have
  
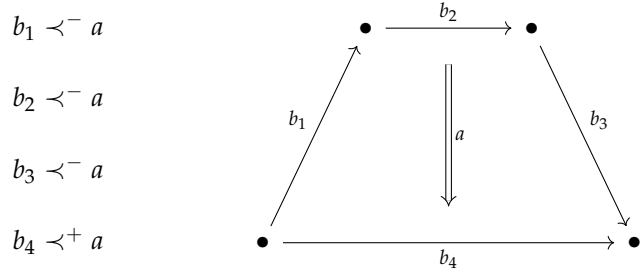
\[\begin{tikzcd}[row sep = small]
    {b_1 \prec^{-} a} &&& \bullet && \bullet \\
    {b_2 \prec^{-} a} \\
    {b_3 \prec^{-} a} \\
    {b_4 \prec^{+} a} && \bullet &&&& \bullet
    \arrow["{b_1}", from=4-3, to=1-4]
    \arrow[""{name=0, anchor=center, inner sep=0}, "{b_2}", from=1-4, to=1-6]
    \arrow["{b_3}", from=1-6, to=4-7]
    \arrow[""{name=1, anchor=center, inner sep=0}, "{b_4}"', from=4-3, to=4-7]
    \arrow["a", shorten <=13pt, shorten >=13pt, Rightarrow, from=0, to=1]
  \end{tikzcd}\]

  \noindent In order to define formally DFCs, we first need the notion of \emph{positive-to-one poset} (POP).

  \begin{defin}[Positive-to-one poset]\label{defin:OGP}
    A positive-to-one poset consists of:
    \begin{itemize}
      \item A finite set of elements $P$.
      \item A gradation $\dim : P \to \N$.
      \item Two binary relations $\prec^{-}$ and $\prec^{+}$ on $P$,
        and we let $x \prec y$ iff $x \prec^{-} y$ or $x \prec^{+} y$.
    \end{itemize}
    With the following properties:
    \begin{itemize}
      \item $\forall x,\,y\in P,\quad y \prec x \to \dim(x) = \dim(y)+1$.
      \item $\forall x,\,y\in P,\quad  \neg\,(y \prec^{-} x \,\wedge\, y \prec^{+} x)$.
      \item $\forall x \in P,\quad  \dim(x)\ge1 \rightarrow (\exists!y,\,y\prec^{+}x) \wedge (\exists y,\,y\prec^{-}x)$.
    \end{itemize}
    In particular: $\prec$, $\prec^{-}$ and $\prec^{+}$ are asymmetric, and the reflexive transitive closure of $\prec$ equips $P$ with a structure of partially ordered set, such that $\dim$ is an increasing map.\\
    Following the conventions of \cite{zawadowski2023positive}, for $x\in P$, we denote 
    $$\delta(x) := \{y\in P \bigm\vert y \prec^{-} x\}$$
    and when $\dim(x) \ge 1$,
    $$\gamma(x) := \{y\in P \bigm\vert y \prec^{+} x\}$$
    because of the third property, $\gamma(x)$ is always a singleton, hence we sometimes identify $\gamma(x)$ with its unique element, which we call the \emph{target} of $x$.
    For $k\in\N$, we also denote $$P_k := \dim^{-1}\left(\{k\}\right),\,\, P_{\ge k} := \bigcup_{i\ge k}P_i,\,\, P_{>k} := \bigcup_{i>k}P_i$$ and we let $\dim(P) := \max\{\dim(x)\}_{x\in P}$ be the \emph{dimension} of $P$.
  \end{defin}

  \para From now on, we will use a kind of \textsc{Hasse} diagrams to represent POPs.
  The convention will be as follows:
  
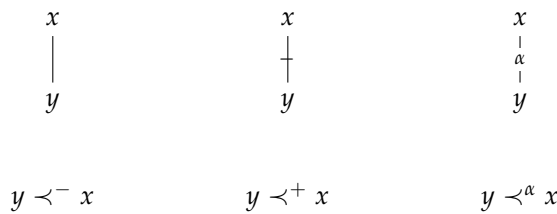
\[\begin{tikzcd}
    x && x && x \\
    y && y && y \\
    {y\prec^{-}x} && {y\prec^{+}x} && {y\prec^{\alpha}x}
    \arrow[no head, from=2-1, to=1-1]
    \arrow["\shortmid"{marking}, no head, from=1-3, to=2-3]
    \arrow["\alpha"{description}, no head, from=1-5, to=2-5]
  \end{tikzcd}\]

  \begin{defin}[Dendritic face complex]\label{defin:DFC}
    A dendritic face complex is a positive-to-one poset $C$, satisfying the following extra axioms:
    \begin{itemize}
      \item (\emph{greatest element})\\
        There is a greatest element in $C$, for the partial order induced by $\prec$.
      \item (\emph{oriented thinness})\\
        For $z \prec y \prec x$ in $P$, there is a unique $y'\neq y$ in $P$ such that $z \prec y' \prec x$.
        Hence there is a lozenge as in \ref{fig:dfc_def_lozenge} below. Moreover, we ask for the \emph{sign rule} $\alpha \beta = - \alpha' \beta'$ to be satisfied.
        When finding such a $y'$ we say that we \emph{complete the half lozenge} $z \prec y \prec x$.
      \item (\emph{acyclicity})\\
        For $x\in P_1$, $\delta(x)$ is a singleton.\\
        Let $x\in P_{\ge1}$, then $\delta(x)\neq\emptyset$ and there is no cycle as in \ref{fig:dfc_def_cycle} below.
    \end{itemize}
    \begin{figure}[H]
      \begin{minipage}{.5\textwidth}
        \[\begin{tikzcd}[row sep = 15pt]
          & x \\
          y && {y'} \\
          & z
          \arrow["\alpha"{description}, no head, from=1-2, to=2-1]
          \arrow["\beta"{description}, no head, from=2-1, to=3-2]
          \arrow["{\alpha'}"{description}, no head, from=1-2, to=2-3]
          \arrow["{\beta'}"{description}, no head, from=2-3, to=3-2]
        \end{tikzcd}\]
        \captionof{figure}{Lozenge}
        \label{fig:dfc_def_lozenge}
      \end{minipage}
      \begin{minipage}{.5\textwidth}
        \[\begin{tikzcd}[row sep = 15pt, column sep = small]
          && x \\
          {y = y_1} & {y_2} & {y_3} & {y_p} & {y_1} \\
          {\gamma(y_2)} & {\gamma(y_3)} & {\gamma(y_4)} & {\gamma(y_1)}
          \arrow[no head, from=2-1, to=3-1]
          \arrow["\shortmid"{marking}, no head, from=2-2, to=3-1]
          \arrow[no head, from=2-2, to=3-2]
          \arrow["\shortmid"{marking}, no head, from=3-2, to=2-3]
          \arrow[no head, from=2-3, to=3-3]
          \arrow["\shortmid"{marking}, dotted, no head, from=2-4, to=3-3]
          \arrow[no head, from=2-4, to=3-4]
          \arrow["\shortmid"{marking}, no head, from=3-4, to=2-5]
          \arrow[curve={height=6pt}, no head, from=1-3, to=2-1]
          \arrow[no head, from=1-3, to=2-2]
          \arrow[no head, from=1-3, to=2-3]
          \arrow[no head, from=1-3, to=2-4]
          \arrow[curve={height=-6pt}, no head, from=1-3, to=2-5]
        \end{tikzcd}\]
        \captionof{figure}{Cycle}
        \label{fig:dfc_def_cycle}
      \end{minipage}
    \end{figure}
  \end{defin}

  \begin{prop}[Tree structure]\label{prop:DFC_tree_structure}
    Let $P$ be a DFC and $x\in P_{\ge1}$. There is a rooted tree structure on $\delta(x)$, given by the following data:
    \begin{itemize}
      \item The set of nodes is $\delta(x)$.
      \item For each $y\in\delta(x)$, $A(y) := \delta(y)$.
      \item There is a triplet $y \branch{z} y'$ iff there is a lozenge as in \ref{fig:tree_structure_triplet} below.
      \item When $\dim(x)\ge 2$, the root is given, as in \cite{zawadowski2023positive}, by the unique $\rho(x)$ completing the lozenge as in \ref{fig:tree_structure_root} below.
    \end{itemize}
    \begin{figure}[H]
      \begin{minipage}{.5\textwidth}
        \[\begin{tikzcd}[column sep = small, row sep = 15pt]
          & x \\
          y && {y'} \\
          & z
          \arrow[no head, from=1-2, to=2-1]
          \arrow[no head, from=2-1, to=3-2]
          \arrow[no head, from=1-2, to=2-3]
          \arrow["\shortmid"{marking}, no head, from=2-3, to=3-2]
        \end{tikzcd}\]
        \captionof{figure}{Triplet}
        \label{fig:tree_structure_triplet}
      \end{minipage}
      \begin{minipage}{.5\textwidth}
        \[\begin{tikzcd}[column sep = small, row sep = 15pt]
          & x \\
          {\rho(x)} && {\gamma(x)} \\
          & {\gamma^2(x)}
          \arrow[no head, from=1-2, to=2-1]
          \arrow["\shortmid"{marking}, no head, from=2-1, to=3-2]
          \arrow["\shortmid"{marking}, no head, from=1-2, to=2-3]
          \arrow["\shortmid"{marking}, no head, from=2-3, to=3-2]
        \end{tikzcd}\]
        \captionof{figure}{Root}
        \label{fig:tree_structure_root}
      \end{minipage}
    \end{figure}
  \end{prop}
  \begin{proof}
    If $\dim(x) = 1$, then $\delta(x)$ is a singleton, hence there is a structure of tree as stated.\\
    Suppose that $\dim(x)>1$.
    By \hyperref[defin:DFC]{oriented thinness}, we may define the root of $\delta(x)$ as above.\\
    If there is a triplet $y \branch{z} y'$ in $\delta(x)$, then $y'$ is uniquely determined by the uniqueness of lozenge completion.\\
    Let $y\in \delta(x)$, then by successively completing lozenges (from left to right in the following diagram) we obtain a configuration as below:
    
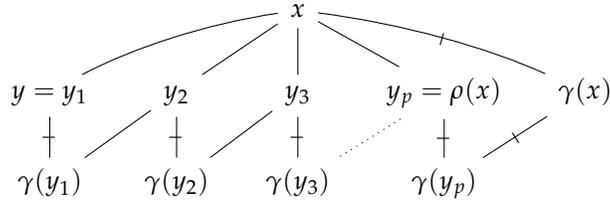
\[\begin{tikzcd}[column sep = small, row sep = 15pt]
      && x \\
      {y = y_1} & {y_2} & {y_3} & {y_{p} = \rho(x)} & {\gamma(x)} \\
      {\gamma(y_1)} & {\gamma(y_2)} & {\gamma(y_3)} & {\gamma(y_{p})}
      \arrow["\shortmid"{marking}, no head, from=2-1, to=3-1]
      \arrow[no head, from=2-2, to=3-1]
      \arrow["\shortmid"{marking}, no head, from=2-2, to=3-2]
      \arrow[no head, from=3-2, to=2-3]
      \arrow["\shortmid"{marking}, no head, from=2-3, to=3-3]
      \arrow[dotted, no head, from=2-4, to=3-3]
      \arrow[curve={height=6pt}, no head, from=1-3, to=2-1]
      \arrow[no head, from=1-3, to=2-2]
      \arrow[no head, from=1-3, to=2-3]
      \arrow[no head, from=1-3, to=2-4]
      \arrow["\shortmid"{marking}, no head, from=2-4, to=3-4]
      \arrow["\shortmid"{marking}, no head, from=3-4, to=2-5]
      \arrow["\shortmid"{marking}, curve={height=-6pt}, no head, from=1-3, to=2-5]
    \end{tikzcd}\]
    That is, we keep completing lozenges while the completion of $\gamma(y_i) \prec^{+} y_i \prec^{-} x$ yields $\gamma(y_i) \prec^{-} y_{i+1} \prec^{-} x$.
    Because of \hyperref[defin:DFC]{acyclicity}, there is some $p\ge 1$ such that the completion of $\gamma(y_p) \prec^{+} y_p \prec^{-} x$ must yield the rightmost lozenge as above.
    Hence, this gives a descending path from $y$ to the root of $\delta(x)$.
    Because of uniqueness of lozenge completion, this is the unique such path.
  \end{proof}

  \begin{defin}[morphism of DFC]\label{defin:DFC_hom}
    Let $C$ and $D$ be two dendritic faces structures.
    Then a morphism $f : C \to D$ of DFCs corresponds to the data of such a map $f$ between the underlying sets of $C$ and $D$, such that:
    \begin{itemize}
      \item $f$ preserves the dimension.
      \item When $y \prec^{+} x$ in $C$, $f(y)\prec^{+}f(x)$ in $D$.
      \item When $y \prec^{-} x$ in $C$, $f(y)\prec^{-}f(x)$ in $D$.
      \item $f_x := f\vert_{\delta(x)} : \delta(x) \to \delta(f(x))$ is a bijection.
    \end{itemize}
  \end{defin}

  \begin{rmk}\label{rmq:DFC_hom_is_iso_on_image}
    It follows from the definition that every morphism $f : C \to D$ induces an isomorphism onto its image $\cl{f(\top)}$ where $\top$ denotes the greatest element of $C$, and $\mathsf{cl}$ the downward closure.
  \end{rmk}

  \begin{defin}\label{defin:OGP_Lambda_Gamma}
    Let $C$ be an POP of dimension $n$. For $k\in\intInter{0}{n}$, we introduce the two following sets:
    $$\Gamma_k := \gamma(C_{k+1}) \qquad \mbox{and its complement} \qquad \Lambda_k := C_k \setminus \gamma(C_{k+1})$$
  \end{defin}

  
\section{Positive opetopes}
\label{sec:ZPO}

  In this section, we describe positive opetopes, as they are defined by \textsc{Zawadowski} in \cite{zawadowski2023positive}.
  They are another formalism dedicated to encoding opetopes, in the wider context of opetopic cardinals.
  We need first a notion of \emph{positive hypergraphs}.

  \begin{defin}[Positive hypergraph]\label{defin:PH}
    A positive hypergraph consists of:
    \begin{itemize}
      \item A family of finite sets $(S_k)_{k\in\N}$ of \emph{faces} such that for $k$ large enough, $S_k = \emptyset$.
      \item For all $k\in N$, a function $\gamma_k : S_{k+1} \to S_k$ and a total relation $\delta_k : S_{k+1} \to S_k$.
        We moreover ask for $\delta_0$ to be functional.
    \end{itemize}
    We introduce the notations $S_{\ge k}$ for $\bigcup_{i\ge k}S_i$ and $S_{>k}$ for $\bigcup_{i>k}S_i$.
    Following \cite{zawadowski2023positive}, we will omit indices for $\delta_i$'s and $\gamma_i$'s, and denote for any $x\in S$, $\gamma^{(k)}(x)$ the iterate $\gamma^q(x)$ such that $\gamma^{q}(x)\in S_k$.\\
    We also let, for $x\in S$, $\iota(x) := \delta\delta(x) \cap \gamma\delta(x)$.
  \end{defin}

  \begin{defin}[Morphism of positive hypergraphs]\label{defin:PH_hom}
    Let $S = (S_k)_{k\in\N}$ and $T = (T_k)_{k\in\N}$ be two positive hypergraphs, then a \emph{morphism} between $S$ and $T$ corresponds to a family of maps
    $(f_k : S_k \to T_k)_{k\in\N}$ such that:
    \begin{itemize}
      \item $\forall x \in S_{\ge1}$, the restriction $f_x := f\vert_{\delta(x)} : \delta(x) \to \delta(f(x))$ is a bijection.
      \item $\forall x \in S_{\ge1}$, $\gamma(f(x)) = f(\gamma(x))$.
    \end{itemize}
  \end{defin}

  \para The \hyperref[defin:OGP]{positive-to-one posets} are the analogue of positive hypergraphs in the formalism of dendritic face structures.
  More precisely:
  let $\mathbf{POP}$ be the category of positive-to-one posets, and $\mathbf{pHg}$ the category of positive hypergraphs.
  Then we may associate to a POP $\left(P = \coprod_{k\in\N}P_k,\,\prec^{+},\,\prec^{-} \right)$ a positive hypergaph $F(P)$ such that:
  \begin{itemize}
    \item For all $k\in\N$, the set of $k$-dimensional faces of $F(P)$ is defined as $P_k$.
    \item For all $k>0$, we define the function $\gamma_k$ as $\gamma\vert : P_{k+1} \to P_k$.
    \item For all $k>0$, we define the total relation $\delta_k$ as $\delta\vert : P_{k+1} \to \mathcal{P}(P_k)$.
  \end{itemize}
  Conversely, to any positive hypergraph $(S = (S_k)_{k\in\N},\,\gamma,\,\delta)$, we may associate a positive-to-one poset $G(S)$ defined as follows:
  \begin{itemize}
    \item The set of elements is $\coprod_{k\in\N}S_k$.
    \item For all $k$, the gradation $\dim$ sends any $x\in S_k$ on $k$.
    \item $\prec^{-}$ is defined by $y \prec^{-} x$ iff $y \in \delta(x)$, and $\prec^{+}$ by $y \prec^{+} x$ iff $y = \gamma(x)$.
  \end{itemize}
  We may also associate to any morphism $f : P \to Q$ of POPs a morphism of positive hypergaphs $F(f) : F(P) \to F(Q)$ defined by $F(f)_k := f\vert : P_k \to Q_k$.
  Conversely, to any morphism of positive hypergraphs $g : S \to T$ we may associate a morphism of POPs $G(g) : G(S) \to G(T)$ defined by $G(g)\vert_{S_k}^{T_k} := g_k$.
  A straightforward check yields the following result:

  \begin{thm}\label{thm:POP_pHg_equiv}
    The functors $\mathbf{POP} \overunderset{F}{G}{\rightleftarrows} \mathbf{pHg}$ form an equivalence of categories.
  \end{thm}

  \para From now on, we will make implicit use of this correspondance. For instance we may use the notation $\prec^{-}, \prec^{+}$ and $\prec$ for positive hypergraphs.

  \begin{defin}[$\triangleleft^{S_k,\,+},\,\triangleleft^{S_k,\,-},\,<^{S_k,\,+},\,<^{S_k,\,-}$]\label{defin:PFS_path_order}
    We define the following relations:
    \begin{itemize}
      \item $<^{S_0,\,-}$ is the empty relation.\\
        For $k > 0$, $<^{S_k,\,-}$ is the transitive closure of
        the relation $\triangleleft^{S_k,\,-}$ on $S_k$, such that $x \triangleleft^{S_k,-} x'$ iff $\gamma(x) \in \delta(x')$. 
        We write $x \bot^{-} x'$ iff either $x <^{-} x'$ or $x' <^{-} x$, and we write $x \le^{-} x'$ iff either $x = x'$
        or $x <^{-} x'$.\\
        A \emph{lower path} is a sequence $x_0 \succ^{+} y_0 \prec^{-} \cdots \prec^{-} x_{p-1} \succ^{+} y_{p} \prec^{-} x_{p}$.
      \item $<^{S_k,\,+}$ is the transitive closure of the relation $\triangleleft^{S_k,+}$ on $S_k$, such that
        $x \triangleleft^{S_k,+} x'$ iff there is $w \in S_{k+1}$, such that $x \in \delta(w)$ and $\gamma(w) = x'$.
        We write $x \bot^{+} x'$ iff either $x <^{+} x'$ or $x' <^{+} x$, and we write $x \le^{+} x'$ iff either $x = x'$
        or $x <^{+} x'$.\\
        An \emph{upper path} is a sequence $y_0 \prec^{-} x_1 \succ^{+} y_1 \prec^{-} \cdots \prec^{-} x_p \succ^{+} y_p$.
    \end{itemize}
  \end{defin}

  \begin{defin}[Opetopic cardinal, positive opetope]\label{defin:PFS}
    An opetopic cardinal corresponds to the data of a positive hypergraph $S$ satisfying the following axioms:
    \begin{itemize}
      \item (\emph{Globularity})\\
        For $x \in S_{\ge2}$,
        \vspace{-10pt}
        \begin{mathpar}
          \gamma\gamma(x) = \gamma\delta(x) \setminus \delta\delta(x)
          \and
          \delta\gamma(x) = \delta\delta(x) \setminus \gamma\delta(x).
        \end{mathpar}
        \vspace{-20pt}
      \item (\emph{Strictness})\\
        For $k\in\N$, $<^{S_k,\,+}$ is a strict order and $<^{S_0,\,+}$ is linear.
      \item (\emph{Disjointness})\\
        For $k>0$,\quad $\bot^{S_k,\,-} \cap \bot^{S_k,\,+} = \emptyset$.
      \item (\emph{Pencil linearity})\\
        For any $k>0$ and $y\in S_{k-1}$, the sets below are linearly ordered by $<^{S_k,\,+}$.
        \begin{mathpar}
          \{x \in S_k \bigm\vert y = \gamma(x)\} \and \mbox{and} \and \{x \in S_k \bigm\vert y \in \delta(x)\}.
        \end{mathpar}
        \vspace{-20pt}
    \end{itemize}
    An opetopic cardinal $S$ is called \emph{principal} if for all $k\le\dim(S)$, $|S_k\setminus\delta(S_{k+1})| = 1$.\\
    A \emph{positive opetope} is a principal opetopic cardinal.
  \end{defin}

  
\section{From dendritic face complexes to \textsc{Zawadowski}'s positive opetopes}
\label{sec:dfc_to_zpo}

\para Let $\left(C = \coprod_{k\in\N} C_k,\,\prec^{+},\,\prec^{-}\right)$ be a dendritic face complex.
Recall from \refthm{thm:POP_pHg_equiv} that $C$ may be given a structure of positive hypergraph. Denote by $n$ the dimension of its greatest element $\omega$.
Our goal is to prove that $C$ is indeed a positive opetope, in the sense of \textsc{Zawadowski}.
We prove the required properties in the following order: \hyperref[defin:PFS]{globularity}, \hyperref[defin:PFS]{strictness} and \hyperref[defin:PFS]{principality}, \hyperref[defin:PFS]{pencil linearity} and then \hyperref[defin:PFS]{disjointness}.
In what follows, $<^{+},\, <^{-},\, \triangleleft^{+}$ and $\triangleleft^{-}$ refer to the corresponding notations introduced in \refdefin{defin:PFS_path_order}.

\begin{prop}[globularity]\label{prop:DFC_glob}
  The positive hypergraph $C$ satisfies the property of \hyperref[defin:PFS]{globularity}.
  That is:\\[-10pt]
  $$
  \gamma\gamma = \gamma\delta - \delta\delta \qquad\mbox{and}\qquad \delta\gamma = \delta\delta - \gamma\delta
  $$
\end{prop}
\begin{proof}
  We first prove that $\gamma\gamma = \gamma\delta - \delta\delta$:\\
  Let $k\ge2$ and $a\in C_k$, and let $b\in\delta(a)$.
  Then we have $\gamma(b)\prec^{+}b\prec^{-}a$, and hence two possible lozenge completions:
  
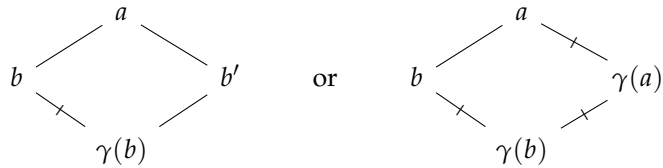
\[\begin{tikzcd}[column sep = scriptsize, row sep = small]
    & a &&&& a \\
    b && {b'} & {\mbox{or}} & b && {\!\!\gamma(a)} \\
    & {\gamma(b)} &&&& {\gamma(b)}
    \arrow[no head, from=1-2, to=2-1]
    \arrow["\shortmid"{marking}, no head, from=2-1, to=3-2]
    \arrow[no head, from=1-2, to=2-3]
    \arrow[no head, from=2-3, to=3-2]
    \arrow[no head, from=1-6, to=2-5]
    \arrow["\shortmid"{marking}, no head, from=2-5, to=3-6]
    \arrow["\shortmid"{marking}, no head, from=1-6, to=2-7]
    \arrow["\shortmid"{marking}, no head, from=2-7, to=3-6]
  \end{tikzcd}\]
  In the first case, $\gamma(b) \in \delta\delta(a)$.
  In the second case, $\gamma(b) = \gamma^2(a)$. Hence $\gamma\delta(a) - \delta\delta(a) \subseteq \gamma\gamma(a)$
  For the converse inclusion, we observe that there is only one possibility of lozenge completion:
  
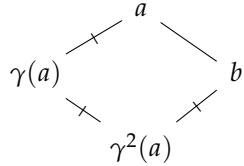
\[\begin{tikzcd}[column sep=small, row sep = small]
    & a \\
    {\gamma(a)\!\!} && b \\
    & {\gamma^2(a)}
    \arrow["\shortmid"{marking}, no head, from=1-2, to=2-1]
    \arrow["\shortmid"{marking}, no head, from=2-1, to=3-2]
    \arrow[no head, from=1-2, to=2-3]
    \arrow["\shortmid"{marking}, no head, from=2-3, to=3-2]
  \end{tikzcd}\]
  Hence, $\gamma^2(a)\in \gamma\delta(a)$ and because of the sign rule, it is impossible to have $\gamma^2(a)\in\delta\delta(a)$.
  So by double inclusion, we have shown that $\gamma\gamma(a) = \gamma\delta(a) - \delta\delta(a)$.

  \noindent We now prove the second equation: $\delta\gamma = \delta\delta - \gamma\delta$.
  Let $k\ge2$, $a\in C_k$, and $c\prec^{-}b\prec^{-}a$.
  Then we have two possible lozenge completions:
  
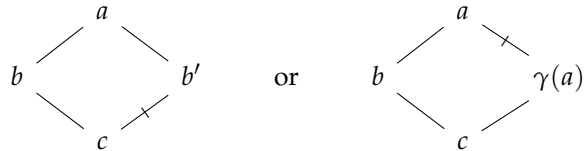
\[\begin{tikzcd}[column sep=scriptsize, row sep = small]
    & a &&&& a \\
    b && {b'} & {\mbox{or}} & b && {\!\!\gamma(a)} \\
    & c &&&& c
    \arrow[no head, from=1-2, to=2-1]
    \arrow[no head, from=2-1, to=3-2]
    \arrow[no head, from=1-2, to=2-3]
    \arrow["\shortmid"{marking}, no head, from=2-3, to=3-2]
    \arrow[no head, from=1-6, to=2-5]
    \arrow[no head, from=2-5, to=3-6]
    \arrow["\shortmid"{marking}, no head, from=1-6, to=2-7]
    \arrow[no head, from=2-7, to=3-6]
  \end{tikzcd}\]
  In the first case, $c\in \gamma\delta(a)$ and in the second one, $c\in\delta\gamma(a)$.
  Hence $\delta\delta(a) - \gamma\delta(a) \subseteq \delta\gamma(a)$.
  On the other hand, if $c\prec^{-}\gamma(a)\prec^{+}a$, then the only possible shape for lozenge completion is
  
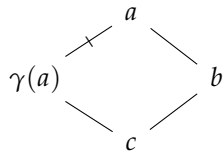
\[\begin{tikzcd}[column sep=scriptsize, row sep = small]
    & a \\
    {\gamma(a)\!\!} && b \\
    & c
    \arrow["\shortmid"{marking}, no head, from=1-2, to=2-1]
    \arrow[no head, from=2-1, to=3-2]
    \arrow[no head, from=1-2, to=2-3]
    \arrow[no head, from=2-3, to=3-2]
  \end{tikzcd}\]
  Hence $c\in\delta\delta(a)$, and because of the sign rule, $c\notin\gamma\delta(a)$. \\
  Whence the second equality: $\delta\gamma(a) = \delta\delta(a) - \gamma\delta(a)$.
\end{proof}

\begin{prop}[strictness, first part]\label{prop:DFC_strictness_1}
  $C$ satisfies the first half of the axiom of \hyperref[defin:PFS]{strictness}.
  That is: $<^{+}$ is a strict partial order.
\end{prop}
\begin{proof}
  We shall prove the following three properties:
  $$
  \begin{array}{ccc}
      \left(k\in\intInter{0}{n}\right) & \mathcal{P}_k : & <^{\scriptscriptstyle C_k,\,+} \mbox{ is a strict order} \\
      \left(k\in\intInter{0}{n-1}\right) & \mathcal{Q}_k : & \forall e\in\Lambda_k,\,\exists!d\in\Lambda_{k+1} \mbox{  s.t.  } e \prec^{-} d \\
      \left(k\in\intInter{0}{n-1}\right) & \mathcal{R}_k : & \forall e\in\Lambda_k,\, e \prec^{-} \gamma^{(k+1)}\omega
  \end{array}
  $$
  by induction on the codimension $n-k$.

  \begin{itemize}
    \item $\boxed{k = n}$ \\
      Since $\omega$ is a greatest element, $C_n = \{\omega\}$. Hence $\mathcal{P}_n$ is clear.
    \item $\boxed{k = n-1}$. \\
      Because $\omega$ is a greatest element, when $a\trigup{n-1}a'$, the only possible situation is the following:
      
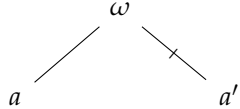
\[\begin{tikzcd}
        & \omega \\
        a && {a'}
        \arrow[no head, from=1-2, to=2-1]
        \arrow["\shortmid"{marking}, no head, from=1-2, to=2-3]
      \end{tikzcd}\]
      Hence a cycle $a_0\trigup{n}a_1\trigup{n}\cdots\trigup{n}a_0$
      must have the form $a_0\trigup{n}a_0$ with $n=0$ and $a_0 = \gamma(\omega)$.
      But $\gamma(\omega)\nprec^{-}\omega$. So there is no such cycle.
      Whence $\mathcal{P}_{n-1}$.

      Because $\omega\in\Lambda_n$, and because it is a greatest element, $\mathcal{Q}_{n-1}$ and $\mathcal{R}_{n-1}$ are clear.

    \item \boxed{\mbox{Induction}}
      We suppose now $\mathcal{P}_{k+1},\,\mathcal{Q}_{k+1}$ and $\mathcal{R}_{k+1}$.
      First, we prove $\mathcal{Q}_k$.
      \begin{itemize}
        \item \textit{Uniqueness}: Let $e\in C_{k}$ and suppose $e\prec^{-}d$ and $e\prec^{-}d'$, with $d,\,d'\in\Lambda_{k+1}$.
          By using $\mathcal{R}_{k+1}$ we may construct a lozenge as in \ref{fig:lozenge_strictness_1} below:
          \begin{figure}[H]
            \hspace{20pt}
            \begin{minipage}{.5\textwidth}
              \[\begin{tikzcd}[row sep = small]
                & {\gamma^{(k+2)}\omega} \\
                d && {d'} \\
                & e
                \arrow[no head, from=1-2, to=2-1]
                \arrow[no head, from=1-2, to=2-3]
                \arrow[no head, from=2-1, to=3-2]
                \arrow[no head, from=2-3, to=3-2]
              \end{tikzcd}\]
              \captionof{figure}{Degenerated lozenge}
              \label{fig:lozenge_strictness_1}
            \end{minipage}%
            \begin{minipage}{.5\textwidth}
              \[\begin{tikzcd}[row sep = small]
                & {c_1} \\
                {d_1} && {d_2} \\
                & e
                \arrow[no head, from=2-1, to=3-2]
                \arrow["\shortmid"{marking}, no head, from=2-1, to=1-2]
                \arrow[no head, from=1-2, to=2-3]
                \arrow[no head, from=2-3, to=3-2]
              \end{tikzcd}\]
              \captionof{figure}{Lozenge}
              \label{fig:strictness_lozenge}
            \end{minipage}
          \end{figure}
          and because of the sign rule, it must be the case that $d = d'$.
        \item \textit{Existence}.
          Let $e\in\Lambda_k$. Since $e$ must be a codimension-$1$ subface of some $d$, we may choose such a $d =: d_1$.
          And because $e$ is not a target, we must have $e\prec^{-} d_1$.
          If $d$ is in $\Lambda_{k+1}$, then we are done.
          Else, let $c_1$ be such that $d_1 = \gamma(c_1)$. Then we may complete $e\prec^{-}d_1\prec^{+}c_1$ into a lozenge as in \ref{fig:strictness_lozenge} above.
          If $d_2\in\Lambda_{k+1}$, we have finished. Else, we continue with $d_2$, taking $c_2$ such that $d_2\prec^{+}c_2$ \textit{etc.}.
          While iterating this construction, we cannot produce a loop as in \ref{fig:strictness_loop} below:
          \begin{figure}[H]
            \hspace{20pt}
            \begin{minipage}{.5\textwidth}
              \[\begin{tikzcd}[column sep = small]
                {c_1} & {c_2} && {c_q} & {c_1} \\
                {d_1} & {d_2} && {d_q} & {d_1} \\
                && e
                \arrow["\shortmid"{marking}, no head, from=1-1, to=2-1]
                \arrow[no head, from=1-1, to=2-2]
                \arrow["\shortmid"{marking}, no head, from=2-2, to=1-2]
                \arrow["\shortmid"{marking}, no head, from=2-5, to=1-5]
                \arrow["\shortmid"{marking}, no head, from=2-4, to=1-4]
                \arrow[no head, from=1-4, to=2-5]
                \arrow[curve={height=12pt}, no head, from=2-1, to=3-3]
                \arrow[curve={height=6pt}, no head, from=2-2, to=3-3]
                \arrow[curve={height=-6pt}, no head, from=2-4, to=3-3]
                \arrow[curve={height=-12pt}, no head, from=2-5, to=3-3]
                \arrow[dotted, no head, from=1-2, to=2-4]
              \end{tikzcd}\]
              \captionof{figure}{Impossible loop}
              \label{fig:strictness_loop}
            \end{minipage}
            \begin{minipage}{.5\textwidth}
              \[\begin{tikzcd}[column sep = small]
                && {\gamma^{(k+2)}(\omega)} \\
                & {d_1} & {d_2} & {d_q} \\
                {e_0} & {e_1} & {e_2} & {e_q}
                \arrow["\shortmid"{marking}, no head, from=3-2, to=2-2]
                \arrow[no head, from=3-1, to=2-2]
                \arrow[no head, from=3-2, to=2-3]
                \arrow["\shortmid"{marking}, no head, from=3-3, to=2-3]
                \arrow[dotted, no head, from=3-3, to=2-4]
                \arrow["\shortmid"{marking}, no head, from=2-4, to=3-4]
                \arrow[no head, from=2-2, to=1-3]
                \arrow[no head, from=1-3, to=2-3]
                \arrow[dotted, no head, from=1-3, to=2-4]
              \end{tikzcd}\]
              \captionof{figure}{Path}
              \label{fig:strictness_path}
            \end{minipage}
          \end{figure}
          because it would contradict $\mathcal{P}_{k+1}$.
          Hence, each time $d_i\notin\Lambda_{k+1}$, the next element $d_{i+1}$ is a new one, and this construction must finish because the poset is finite.
          So when the construction ends, it produces an element $e\prec^{-}d_q\in\Lambda_{k+1}$. Whence the existence.
      \end{itemize}
      
    \vspace{10pt}
    We then prove $\mathcal{R}_k$.
    Let $e\in\Lambda_{k+1}$, then using $\mathcal{Q}_{k}, \mathcal{R}_{k+1}$, we may find some $d\in\Lambda_{k+1}$ such that $e \prec^{-} d \prec^{-} \gamma^{(k+2)}\omega$.
    Since $e$ is not a source, the only possible lozenge completion for this triple yield $e \prec^{-} \gamma^{(k+1)}\omega \prec^{+} \gamma^{(k+2)}\omega$.
    Whence $\mathcal{R}_k$.

    \vspace{10pt}
    Now we prove $\mathcal{P}_k$.
    The strategy is to replace a related pair $e\trigup{k}e'$ by a path\\
    $e = e_0 \trigup{k} e_1 \trigup{k} \cdots \trigup{k} e_q = e'$ as in \ref{fig:strictness_path} above:
    Indeed, if we are able to do this, every cyclic path $e_0\trigup{k}e_1\trigup{k}\cdots\trigup{k}e_q = e_0$ will induce a longer path
    $e_0 = e'_0\trigup{k}e'_1\trigup{k}\cdots\trigup{k}e'_m = e_0$ as above, which cannot exist because it would imply the existence of a cycle in the tree structure of $\delta(\gamma^{(k+2)}\omega)$.
    First, notice that if we know each $d_i$ is in $\Lambda_{k+1}$ then we know that $d_i \prec^{-} \gamma^{(k+2)}\omega$ (by $\mathcal{R}_{k+1}$).
    So we will be able to conclude if we prove \ref{lem:DFC_path_in_Lambda} below.
      
    \begin{lem}\label{lem:DFC_path_in_Lambda}
      if $e\trigup{k}e'$, then we have a sequence as follows, with $d_1,\,d_2,\,\cdots,\,d_q\in\Lambda_{k+1}$.
      
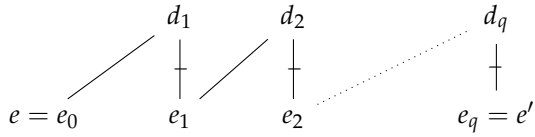
\[\begin{tikzcd}
        & {d_1} & {d_2} && {d_q} \\
        {e = e_0} & {e_1} & {e_2} && {e_q = e'}
        \arrow["\shortmid"{marking}, no head, from=2-2, to=1-2]
        \arrow[no head, from=2-1, to=1-2]
        \arrow[no head, from=2-2, to=1-3]
        \arrow["\shortmid"{marking}, no head, from=2-3, to=1-3]
        \arrow[dotted, no head, from=2-3, to=1-5]
        \arrow["\shortmid"{marking}, no head, from=1-5, to=2-5]
      \end{tikzcd}\]
    \end{lem}
    \textit{Proof.}
      If $e\prec^{-}d\succ^{+}e'$ with $d\in\Lambda_{k+1}$, then the result is proven.
      So we may suppose that $d$ is the target of some element: $d = \gamma(c)$, as in \ref{fig:dfc_to_zpo_hasse_1} below.
      \begin{figure}[H]
        \begin{minipage}{.5\textwidth}
          \[\begin{tikzcd}[row sep = small]
            & c \\
            & d \\
            e && {e'}
            \arrow[no head, from=3-1, to=2-2]
            \arrow["\shortmid"{marking}, no head, from=2-2, to=3-3]
            \arrow["\shortmid"{marking}, no head, from=2-2, to=1-2]
          \end{tikzcd}\]
          \captionof{figure}{}
          \label{fig:dfc_to_zpo_hasse_1}
        \end{minipage}
        \begin{minipage}{.5\textwidth}
          \[\begin{tikzcd}[row sep = small]
            & c \\
            d & {d_1} \\
            e
            \arrow[no head, from=3-1, to=2-2]
            \arrow[no head, from=2-2, to=1-2]
            \arrow["\shortmid"{marking}, no head, from=1-2, to=2-1]
            \arrow[no head, from=2-1, to=3-1]
          \end{tikzcd}\]
          \captionof{figure}{}
          \label{fig:dfc_to_zpo_hasse_2}
        \end{minipage}
      \end{figure}
      First, complete the left half lozenge with some $d_1$ as in \ref{fig:dfc_to_zpo_hasse_2} above.
      
      Then we consider its target $\gamma(d_1)$, and complete $\gamma(d_1)\prec^{+}d_1\prec^{-}c$ as a lozenge. We have two possibilities.
      Either the completion is of the type 1 as in \ref{fig:dfc_to_zpo_hasse_3} below, or it of the type 2 as in \ref{fig:dfc_to_zpo_hasse_4} below:
      \begin{figure}[H]
        \begin{minipage}{.5\textwidth}
          \[\begin{tikzcd}[row sep = small]
            & c \\
            d & {d_1} & {d_2 = d} \\
            e & {\gamma(d_1) = e'}
            \arrow[no head, from=3-1, to=2-2]
            \arrow[no head, from=2-2, to=1-2]
            \arrow["\shortmid"{marking}, no head, from=1-2, to=2-1]
            \arrow[no head, from=2-1, to=3-1]
            \arrow["\shortmid"{marking}, no head, from=2-2, to=3-2]
            \arrow["\shortmid"{marking}, no head, from=1-2, to=2-3]
            \arrow["\shortmid"{marking}, no head, from=2-3, to=3-2]
          \end{tikzcd}\]
          \captionof{figure}{Type 1}
          \label{fig:dfc_to_zpo_hasse_3}
        \end{minipage}
        \begin{minipage}{.5\textwidth}
          \[\begin{tikzcd}[row sep = small]
            & c \\
            d & {d_1} & {d_2} \\
            e & {\gamma(d_1)}
            \arrow[no head, from=3-1, to=2-2]
            \arrow[no head, from=2-2, to=1-2]
            \arrow["\shortmid"{marking}, no head, from=1-2, to=2-1]
            \arrow[no head, from=2-1, to=3-1]
            \arrow["\shortmid"{marking}, no head, from=2-2, to=3-2]
            \arrow[no head, from=1-2, to=2-3]
            \arrow[no head, from=2-3, to=3-2]
          \end{tikzcd}\]
          \captionof{figure}{Type 2}
          \label{fig:dfc_to_zpo_hasse_4}
        \end{minipage}
      \end{figure}
      In the second case, we keep completing lozenges on the right until ending with a diagram of the shape of \ref{fig:dfc_to_zpo_hasse_5} below:
      \begin{figure}[H]
        \begin{minipage}{.5\textwidth}
          \[\begin{tikzcd}[column sep = small]
            && c \\
            d & {d_1} & {d_2} & {d_q} & d \\
            e & {\gamma(d_1)} & {\gamma(d_2)} & {\gamma(d_q) = e'}
            \arrow[no head, from=3-1, to=2-2]
            \arrow[no head, from=2-1, to=3-1]
            \arrow["\shortmid"{marking}, no head, from=2-2, to=3-2]
            \arrow[no head, from=2-3, to=3-2]
            \arrow["\shortmid"{marking}, no head, from=2-3, to=3-3]
            \arrow[dotted, no head, from=3-3, to=2-4]
            \arrow["\shortmid"{marking}, no head, from=2-4, to=3-4]
            \arrow["\shortmid"{marking}, no head, from=2-5, to=3-4]
            \arrow["\shortmid"{marking}, curve={height=12pt}, no head, from=1-3, to=2-1]
            \arrow[curve={height=-6pt}, no head, from=2-2, to=1-3]
            \arrow[no head, from=1-3, to=2-3]
            \arrow[curve={height=-6pt}, no head, from=1-3, to=2-4]
            \arrow["\shortmid"{marking}, curve={height=-12pt}, no head, from=1-3, to=2-5]
          \end{tikzcd}\]
          \captionof{figure}{Path}
          \label{fig:dfc_to_zpo_hasse_5}
        \end{minipage}
        \begin{minipage}{.5\textwidth}
          \[\begin{tikzcd}[column sep = small]
            & c \\
            {d_j} & {d_{j+1}} & {d_{i+1} = d_j} \\
            {\gamma(d_j)} & {\gamma(d_{j+1})} & {\gamma(d_{i+1}) = \gamma(d_j)}
            \arrow["\shortmid"{marking}, no head, from=2-1, to=3-1]
            \arrow[no head, from=2-2, to=3-1]
            \arrow["\shortmid"{marking}, no head, from=2-2, to=3-2]
            \arrow[dotted, no head, from=3-2, to=2-3]
            \arrow["\shortmid"{marking}, no head, from=2-3, to=3-3]
            \arrow[no head, from=2-1, to=1-2]
            \arrow[no head, from=1-2, to=2-2]
            \arrow[no head, from=1-2, to=2-3]
          \end{tikzcd}\]
          \captionof{figure}{Cycle}
          \label{fig:dfc_to_zpo_hasse_6}
        \end{minipage}
      \end{figure}
      First, we will never encounter a situation where $d_{i+1} = d_j$ for a $j\leq i$ as in \ref{fig:dfc_to_zpo_hasse_6} above, because it would imply having a cycle
      $$d_j\branch{\gamma(d_j)}d_{j+1}\branch{\gamma(d_{j+1})}\cdots\branch{\gamma(d_i)} d_{i+1} = d_j$$
      hence a cycle in the tree structure of $\delta(c)$.
      And this construction must finish as above because of the finitness of the poset $C$.
      Now, after renaming the elements, we end up with this situation.
      
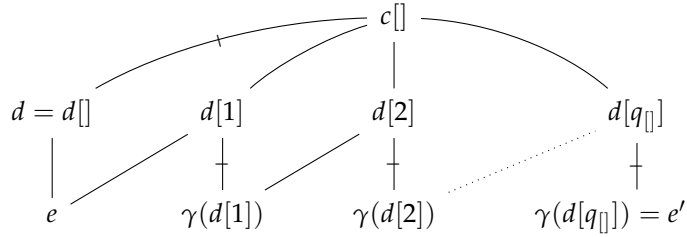
\[\begin{tikzcd}
        && {c[]} \\
        {d = d[]} & {d[1]} & {d[2]} && {d[q_{[]}]} \\
        e & {\gamma(d[1])} & {\gamma(d[2])} && {\hspace{-20pt} \gamma(d[q_{[]}]) = e'}
        \arrow[no head, from=3-1, to=2-2]
        \arrow[no head, from=2-1, to=3-1]
        \arrow["\shortmid"{marking}, no head, from=2-2, to=3-2]
        \arrow[no head, from=2-3, to=3-2]
        \arrow["\shortmid"{marking}, no head, from=2-3, to=3-3]
        \arrow[dotted, no head, from=3-3, to=2-5]
        \arrow["\shortmid"{marking}, no head, from=2-5, to=3-5]
        \arrow["\shortmid"{marking}, curve={height=12pt}, no head, from=1-3, to=2-1]
        \arrow[curve={height=-6pt}, no head, from=2-2, to=1-3]
        \arrow[no head, from=1-3, to=2-3]
        \arrow[curve={height=-12pt}, no head, from=1-3, to=2-5]
      \end{tikzcd}\]
      If each $d[i]$ is in $\Lambda_{k+1}$ then we are done.
      In the other case, for example if $d[2]\prec^{+}c[2]$, then we can once again \textit{unfold}
      
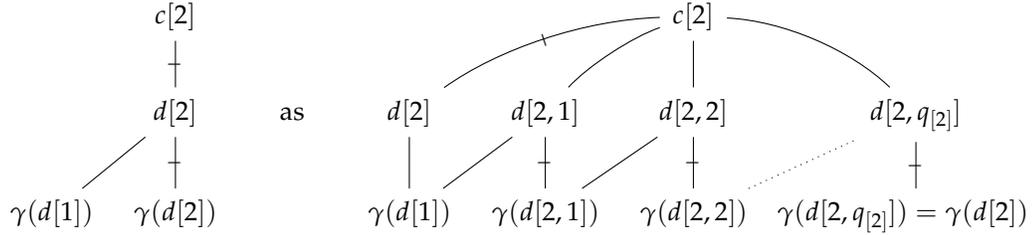
\[\begin{tikzcd}[column sep = tiny]
        & {c[2]} &&&&&& {c[2]} \\
        & {d[2]} && {\mbox{ as }} && {d[2]} & {d[2,1]} & {d[2,2]} && {d[2,q_{[2]}]} \\
        {\gamma(d[1])} & {\gamma(d[2])} &&&& {\gamma(d[1])} & {\gamma(d[2,1])} & {\gamma(d[2,2])} && {\hspace{-10pt} \gamma(d[2,q_{[2]}]) = \gamma(d[2])}
        \arrow[no head, from=3-6, to=2-7]
        \arrow[no head, from=2-6, to=3-6]
        \arrow["\shortmid"{marking}, no head, from=2-7, to=3-7]
        \arrow[no head, from=2-8, to=3-7]
        \arrow["\shortmid"{marking}, no head, from=2-8, to=3-8]
        \arrow[dotted, no head, from=3-8, to=2-10]
        \arrow["\shortmid"{marking}, no head, from=2-10, to=3-10]
        \arrow["\shortmid"{marking}, curve={height=12pt}, no head, from=1-8, to=2-6]
        \arrow[curve={height=-6pt}, no head, from=2-7, to=1-8]
        \arrow[no head, from=1-8, to=2-8]
        \arrow[curve={height=-12pt}, no head, from=1-8, to=2-10]
        \arrow["\shortmid"{marking}, no head, from=1-2, to=2-2]
        \arrow[no head, from=2-2, to=3-1]
        \arrow["\shortmid"{marking}, no head, from=2-2, to=3-2]
      \end{tikzcd}\]
      Now we may iterate this unfolding process. This will produce a tree shaped collection of $d$'s and $c$'s, with relations
      $$d[a_1,\,\cdots,\,a_p]\prec^{+}c[a_1,\,\cdots,\,a_p] \qquad \mbox{and} \qquad c[a_1,\,\cdots,\,a_p]\succ^{-}d[a_1,\,\cdots,\,a_p,\,a]$$ for each $a\in\left[\hspace{-3pt}\left[1,\,q_{[a_1,\,\cdots,\,a_p]}\right]\hspace{-3pt}\right]$.
      Each branch $$d[]\prec^{+}c[]\succ^{-}d[a_1]\prec^{+}c[a_1]\succ^{-}d[a_1,\,a_2]\prec^{+}\cdots\succ^{-}d[a_1,\,\cdots,\,a_p]\prec^{+}\cdots$$
      must be finite, otherwise it would contradict $\mathcal{P}_{k+1}$.
      \textit{i.e.} at some point, all $d[a_1,\,\cdots,\,a_p]$ are in $\Lambda_{k+1}$.
      We thus obtain a path of the desired form, which completes the proof of \reflem{lem:DFC_path_in_Lambda}.
    \qedhere
  \end{itemize}
\end{proof}

\noindent Notice that through the proof above, we have seen the following result (\textit{c.f.} $\mathcal{R}_k$):

\begin{lem}\label{lem:DFC_Lambda_sources_of_iterated_target}
  $\forall k\in\intInter{0}{n-1},\quad \Lambda_k \subseteq \delta(\gamma^{(k+1)}\omega)$.
\end{lem}

\begin{lem}\label{lem:DFC_not_a_source_is_iterated_target}
  Suppose $d\in C_{k}$ is not a source (\textit{i.e.} there is no $c\in C_{k-1}$ such that $d\prec^{-}c$).
  Then $d = \gamma^{(k)}\omega$.
\end{lem}
\begin{proof}
  We will proceed by induction on the codimension of $d$.
  \begin{itemize}
    \item \boxed{k=n} \\
      In this case, $d = \omega$ is the only $n$-dimensional cell.
    \item \boxed{\mbox{Heredity}} \\
      Suppose that we know the result for all cells of dimension $k+1$, and let $d\in C_k$ not being a source.
      Because $\omega$ is a greatest element, we know that there is a chain from $d$ to $\omega$.
      Hence, there is $c_0\in C_{k+1}$ such that $d\prec^{+} c_0$, $i.e.$ $d = \gamma(c)$.
      If $c_0$ is not a source, then we are done by induction hypothesis.
      In the other case, there is $b_0\succ^{-}c$, thus there is a lozenge completion as in \ref{fig:dfc_to_zpo_hasse_7} below:
      \begin{figure}[H]
        \begin{minipage}{.5\textwidth}
          \[\begin{tikzcd}
            & b_0 \\
            c_0 && {c_1} \\
            & d
            \arrow[no head, from=1-2, to=2-1]
            \arrow["\shortmid"{marking}, no head, from=2-1, to=3-2]
            \arrow["\shortmid"{marking}, no head, from=2-3, to=3-2]
            \arrow["\shortmid"{marking}, no head, from=1-2, to=2-3]
          \end{tikzcd}\]
          \captionof{figure}{Lozenge completion}
          \label{fig:dfc_to_zpo_hasse_7}
        \end{minipage}
        \begin{minipage}{.5\textwidth}
          \[\begin{tikzcd}
            {b_0} & {b_1} & {b_2} & {b_{p-1}} \\
            {c_0} & {c_1} & {c_2} & {c_{p-1}} & {c_p} \\
            && d
            \arrow["\shortmid"{marking}, no head, from=1-1, to=2-2]
            \arrow[no head, from=1-1, to=2-1]
            \arrow[no head, from=1-2, to=2-2]
            \arrow[no head, from=1-3, to=2-3]
            \arrow[no head, dotted, from=1-3, to=2-4]
            \arrow["\shortmid"{marking}, no head, from=1-2, to=2-3]
            \arrow[no head, from=2-4, to=1-4]
            \arrow["\shortmid"{marking}, no head, from=2-5, to=1-4]
            \arrow["\shortmid"{marking}, curve={height=12pt}, no head, from=2-1, to=3-3]
            \arrow["\shortmid"{marking}, curve={height=6pt}, no head, from=2-2, to=3-3]
            \arrow["\shortmid"{marking}, no head, from=2-3, to=3-3]
            \arrow["\shortmid"{marking}, curve={height=-6pt}, no head, from=2-4, to=3-3]
            \arrow["\shortmid"{marking}, curve={height=-12pt}, no head, from=2-5, to=3-3]
          \end{tikzcd}\]
          \captionof{figure}{Path}
          \label{fig:dfc_to_zpo_hasse_8}
        \end{minipage}
      \end{figure}
      with $d\prec^{+}c_1$, because $d$ is not a source.
      If $c_1$ is not a source, then we are done by induction hypothesis.
      Else, there is $b_1\in C_{k+1}$, such that $b_1 \succ^{-} c_1$. Hence we may keep completing lozenges as in \ref{fig:dfc_to_zpo_hasse_8} above, until coming accross some $c_p$ which is not a source.
      This proccess must end because of \hyperref[prop:DFC_strictness_1]{strictness} (first part). By induction hypothesis, $c_p$ is an iterated target, and so is $d$.
      \qedhere
  \end{itemize}
\end{proof}

\begin{lem}\label{lem:DFC_source_is_source_of_Lambda}
  If $d$ is a source, then $\exists!c\in\Lambda_{k+1}$ s.t. $d\prec^{-}c$.
\end{lem}
\begin{proof}
  Let $c_0$ be such that $d\prec^{-}c_0$.
  Suppose $c_0\in\Lambda_{k+1}$, then we are done taking $c = c_0$.
  Else as in the proof of \reflem{lem:DFC_not_a_source_is_iterated_target} we may complete lozenges from left to right (refer to \ref{fig:dfc_to_zpo_hasse_9} below),
  until coming accross some $c\in\Lambda_{k+1}$.
  Notice that it must finish because $<^{+}$ is a strict order.
  The uniqueness follows also from the same argument as in the proof of \hyperref[prop:DFC_strictness_1]{$\mathcal{Q}_k$}.
\end{proof}

\begin{lem}\label{lem:DFC_target_is_target_of_Lambda}
  If $d\in C_k$ is a target, then it is the target of a unique $c\in\Lambda_{k+1}$.
\end{lem}
\begin{proof}
  The proof of existence goes by the same kind of construction as above, by filling lozenges as below from left to right until coming accross some $c\in\Lambda_{k+1}$ (refer to \ref{fig:dfc_to_zpo_hasse_10} below).
  The uniqueness comes from the same argument as for the uniqueness in \hyperref[prop:DFC_strictness_1]{$\mathcal{Q}_k$}.
\end{proof}
\begin{figure}[H]
  \begin{minipage}{.5\textwidth}
    \[\begin{tikzcd}
      {b_0} & {b_1} & {b_2} & {b_q} \\
      {c_0} & {c_1} & {c_2} & {c_q} & {c} \\
      && d
      \arrow["\shortmid"{marking}, no head, from=1-1, to=2-1]
      \arrow["\shortmid"{marking}, no head, from=1-2, to=2-2]
      \arrow["\shortmid"{marking}, no head, from=1-3, to=2-3]
      \arrow["\shortmid"{marking}, no head, from=1-4, to=2-4]
      \arrow[no head, from=1-4, to=2-5]
      \arrow[dotted, no head, from=1-3, to=2-4]
      \arrow[no head, from=1-2, to=2-3]
      \arrow[no head, from=1-1, to=2-2]
      \arrow[no head, from=2-5, to=1-4]
      \arrow[curve={height=12pt}, no head, from=2-1, to=3-3]
      \arrow[curve={height=6pt}, no head, from=2-2, to=3-3]
      \arrow[no head, from=2-3, to=3-3]
      \arrow[curve={height=-6pt}, no head, from=2-4, to=3-3]
      \arrow[curve={height=-12pt}, no head, from=2-5, to=3-3]
    \end{tikzcd}\]
    \captionof{figure}{}
    \label{fig:dfc_to_zpo_hasse_9}
  \end{minipage}
  \begin{minipage}{.5\textwidth}
    \[\begin{tikzcd}
      {b_0} & {b_1} & {b_2} & {b_q} \\
      {c_0} & {c_1} & {c_2} & {c_q} & c \\
      && d
      \arrow["\shortmid"{marking}, no head, from=1-1, to=2-1]
      \arrow["\shortmid"{marking}, no head, from=1-2, to=2-2]
      \arrow["\shortmid"{marking}, no head, from=1-3, to=2-3]
      \arrow["\shortmid"{marking}, no head, from=1-4, to=2-4]
      \arrow[no head, from=1-4, to=2-5]
      \arrow[dotted, no head, from=1-3, to=2-4]
      \arrow[no head, from=1-2, to=2-3]
      \arrow[no head, from=1-1, to=2-2]
      \arrow[no head, from=2-5, to=1-4]
      \arrow["\shortmid"{marking}, curve={height=12pt}, no head, from=2-1, to=3-3]
      \arrow["\shortmid"{marking}, curve={height=6pt}, no head, from=2-2, to=3-3]
      \arrow["\shortmid"{marking}, no head, from=2-3, to=3-3]
      \arrow["\shortmid"{marking}, curve={height=-6pt}, no head, from=2-4, to=3-3]
      \arrow["\shortmid"{marking}, curve={height=-12pt}, no head, from=2-5, to=3-3]
    \end{tikzcd}\]
    \captionof{figure}{}
    \label{fig:dfc_to_zpo_hasse_10}
  \end{minipage}
\end{figure}

We mention that \reflem{lem:DFC_Lambda_sources_of_iterated_target} has a converse, although we will not use it later on.

\begin{lem}\label{lem:DFC_sources_of_iterated_target_are_in_Lambda}
  If $d\prec^{-}\gamma^{(k+1)}\omega$, then $d\in\Lambda_k$.
\end{lem}
\begin{proof}
  Suppose $d$ is a target, then because of \reflem{lem:DFC_target_is_target_of_Lambda}, there is $c\in\Lambda_{k+1}$ s.t. $d = \gamma(c)$.
  But because $c$ is in $\Lambda_{k+1}$, we have seen in the proof of \hyperref[prop:DFC_strictness_1]{strictness} that $c\prec^{-}\gamma^{(k+2)}\omega$.
  Hence we are in the situation of \ref{fig:dfc_to_zpo_hasse_11} below,
  which is prohibited by the sign rule, whence $d\in\Lambda_{k+1}$.
\end{proof}

\begin{lem}\label{lem:DFC_iterated_target_is_not_a_source}
  $\gamma^{(k)}\omega$ is not a source. 
\end{lem}
\begin{proof}
  Suppose that it is, then by \reflem{lem:DFC_source_is_source_of_Lambda} it is the source of some $c\in\Lambda_k$.
  Hence we have the lozenge of \ref{fig:dfc_to_zpo_hasse_12} below,
  which is prohibited by the sign rule.
\end{proof}

\begin{figure}[H]
  \begin{minipage}{.5\textwidth}
    \[\begin{tikzcd}[column sep = tiny]
      & {\gamma^{(k+2)}\omega} \\
      {\gamma^{(k+1)}\omega} &&& c \\
      & d
      \arrow["\shortmid"{marking}, no head, from=1-2, to=2-1]
      \arrow[no head, from=2-1, to=3-2]
      \arrow[no head, from=1-2, to=2-4]
      \arrow["\shortmid"{marking}, no head, from=3-2, to=2-4]
    \end{tikzcd}\]
    \captionof{figure}{Forbidden lozenge}
    \label{fig:dfc_to_zpo_hasse_11}
  \end{minipage}
  \begin{minipage}{.5\textwidth}
    \[\begin{tikzcd}[column sep=tiny]
      & {\gamma^{(k+2)}\omega} \\
      {\gamma^{(k+1)}\omega} &&& c \\
      & {\gamma^{(k)}\omega}
      \arrow["\shortmid"{marking}, no head, from=1-2, to=2-1]
      \arrow["\shortmid"{marking}, no head, from=2-1, to=3-2]
      \arrow[no head, from=1-2, to=2-4]
      \arrow[no head, from=2-4, to=3-2]
    \end{tikzcd}\]
    \captionof{figure}{Forbidden lozenge}
    \label{fig:dfc_to_zpo_hasse_12}
  \end{minipage}
\end{figure}

\noindent We also have the following corollary of \reflem{lem:DFC_not_a_source_is_iterated_target}, \reflem{lem:DFC_source_is_source_of_Lambda} and \reflem{lem:DFC_target_is_target_of_Lambda}:
\begin{lem}\label{lem:DFC_source_or_target_of_Lambda}
  If $d\in C_k$ for some $k<n$, it is the source of a unique $c\in\Lambda_{k+1}$ or the target of a unique $c\in\Lambda_{k+1}$.
\end{lem}

\noindent Finally, we may state the following theorem:

\begin{thm}[Sources partition]\label{thm:sources_partition}
  $$(0\le k\le n)\quad C_k\setminus\{\gamma^{(k)}\omega\} = \coprod_{c\in\Lambda_{k+1}}{\delta(c)}$$
\end{thm}
\begin{proof}~\\[-15pt]
  \begin{itemize}
    \item $C_k\setminus\{\gamma^{(k)}\omega\} \subseteq \coprod_{c\in\Lambda_{k+1}}{\delta(c)}$
      is given by \reflem{lem:DFC_not_a_source_is_iterated_target} and \reflem{lem:DFC_source_is_source_of_Lambda}.
    \item $C_k\setminus\{\gamma^{(k)}\omega\} \supseteq \coprod_{c\in\Lambda_{k+1}}{\delta(c)}$
      is given by \reflem{lem:DFC_iterated_target_is_not_a_source}.
  \end{itemize}
\end{proof}

\begin{prop}[principality]\label{prop:DFC_principal}
  The positive hypergraph $C$ is \hyperref[defin:PFS]{principal}.
\end{prop}
\begin{proof}
  It follows directly from \refthm{thm:sources_partition}.
\end{proof}

\begin{prop}[strictness ($<^{C_0,+}$ is total)]\label{prop:DFC_strictness_2}
  $C$ satisfies the second half of the \hyperref[defin:PFS]{strictness} property.
  That is: $<^{C_0,+}$ is a linear total order.
\end{prop}
\begin{proof}
  Let $x\in C_0$.
  From \refthm{thm:sources_partition}, we know that either $x = \gamma^{n}\omega$, or $x$ is the source of a unique $w\in \Lambda_1$.
  By iterating this case disjunction, we may produce a unique path as below:
  
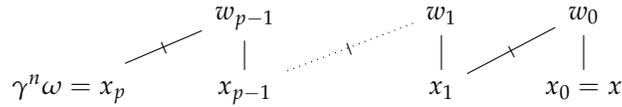
\[\begin{tikzcd}[row sep = small]
    & {w_{p-1}} && {w_1} & {w_0} \\
    {\gamma^n\omega = x_p} & {x_{p-1}} && {x_1} & {x_0 = x}
    \arrow[no head, from=2-5, to=1-5]
    \arrow["\shortmid"{marking}, no head, from=1-5, to=2-4]
    \arrow[no head, from=1-4, to=2-4]
    \arrow["\shortmid"{marking}, dotted, no head, from=1-4, to=2-2]
    \arrow[no head, from=1-2, to=2-2]
    \arrow["\shortmid"{marking}, no head, from=1-2, to=2-1]
  \end{tikzcd}\]
  where the $w_i$'s are in $\Lambda_1$.
  Because of the functionality of $\delta_0$ and \reflem{lem:DFC_target_is_target_of_Lambda}, the path from left to right starting from $\gamma^n\omega$ is also unique (and this is independant of $x$).
  Hence by extending this path as far to the right as possible (which ends because $<_{c_0,+}$ is a strict order), all $0$-dimensional elements must appear somewhere along the path. Which proves that $<^{C_0,+}$ is total.
\end{proof}

\begin{lem}\label{lem:positive_parenthesis_completion}
  For every configuration $e\prec^{\beta}d\prec^{+}c\prec^{-}b$,
  there is a (unique) chain
  $$
  c = c_0 \succ^{-} d_0 \prec^{+} c_1 \succ^{-} d_1 \prec^{+} \cdots \prec^{+} c_p \succ^{-} d_p \prec^{-} c_{p+1} = \gamma(b)
  $$
  with $p\ge0$, as below:
  
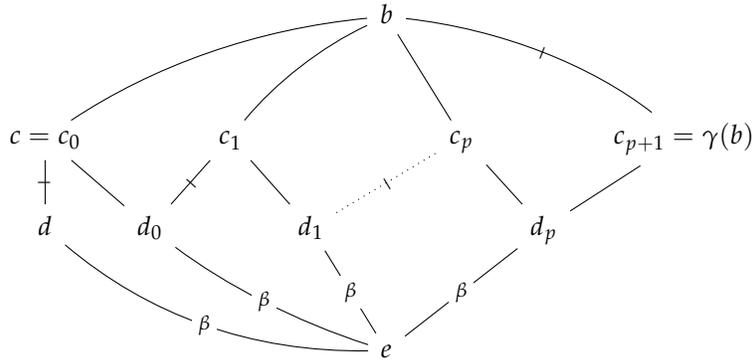
\[\begin{tikzcd}[column sep = small, row sep = 15pt]
    &&&& b \\
    \\
    {c = c_0} && {c_1} &&& {c_p} && {c_{p+1} = \gamma(b)} \\
    d & {d_0} && {d_1} &&& {d_p} \\
    \\
    &&&& e
    \arrow["\shortmid"{marking}, no head, from=3-1, to=4-1]
    \arrow[curve={height=12pt}, no head, from=1-5, to=3-1]
    \arrow[curve={height=6pt}, no head, from=1-5, to=3-3]
    \arrow[no head, from=1-5, to=3-6]
    \arrow["\shortmid"{marking}, curve={height=-12pt}, no head, from=1-5, to=3-8]
    \arrow[no head, from=3-1, to=4-2]
    \arrow["\shortmid"{marking}, no head, from=4-2, to=3-3]
    \arrow[no head, from=3-3, to=4-4]
    \arrow["\shortmid"{marking}, dotted, no head, from=4-4, to=3-6]
    \arrow[no head, from=3-6, to=4-7]
    \arrow[no head, from=4-7, to=3-8]
    \arrow["\beta"{description}, curve={height=18pt}, no head, from=4-1, to=6-5]
    \arrow["\beta"{description}, curve={height=6pt}, no head, from=4-2, to=6-5]
    \arrow["\beta"{description}, no head, from=4-4, to=6-5]
    \arrow["\beta"{description}, no head, from=4-7, to=6-5]
  \end{tikzcd}\]
\end{lem}
\begin{proof}
  By lozenge completion, we may find a unique $d_0$ completing the lozenge $(c,d,e,d_0)$.
  Because of the sign rule, we have $e \prec^{\beta} d_0 \prec^{-} c =: c_0$.
  Then we may find a unique $c_1$ completing the lozenge $(b,c_0,d_0,c_1)$.
  If $d_0 \prec^{-} c_1 \prec^{+} b$, then we are done taking $p=1$.
  Else, we continue this process starting from $d_0$.
  It must finish because of \hyperref[prop:DFC_strictness_1]{strictness}.
\end{proof}

\begin{lem}\label{lem:negative_parenthesis_completion}
  For every configuration $e\prec^{\beta}d\prec^{-}c\prec^{-}b$,
  there is a (unique) chain
  $$
  c = c_0 \succ^{-} d = d_0 \prec^{+} c_1 \succ^{-} d_1 \prec^{+} \cdots \prec^{+} c_{p} \succ^{-} d_{p} \prec^{-} c_{p+1} = \gamma(b)
  $$
  with $p\ge0$ as below:
  
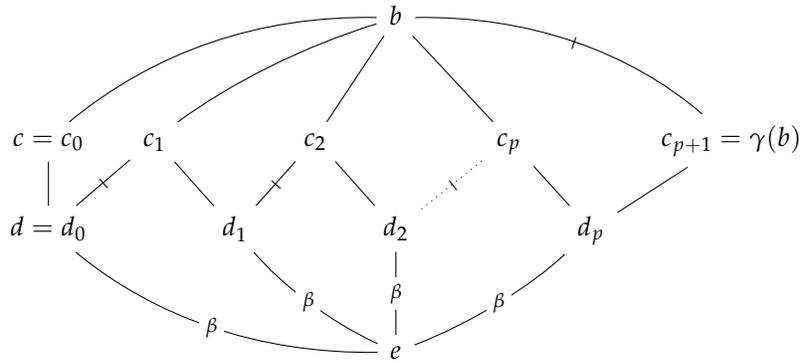
\[\begin{tikzcd}[column sep = small, row sep = 15pt]
    &&&& b \\
    \\
    {c = c_0} & {c_1} && {c_2} &&& {c_p} && {c_{p+1} = \gamma(b)} \\
    {d = d_0} && {d_1} && {d_2} &&& {d_p} \\
    \\
    &&&& e
    \arrow[no head, from=3-1, to=4-1]
    \arrow[curve={height=18pt}, no head, from=1-5, to=3-1]
    \arrow["\beta"{description}, curve={height=18pt}, no head, from=4-1, to=6-5]
    \arrow["\shortmid"{marking}, no head, from=4-1, to=3-2]
    \arrow[no head, from=3-2, to=4-3]
    \arrow["\shortmid"{marking}, no head, from=4-3, to=3-4]
    \arrow[no head, from=3-4, to=4-5]
    \arrow[no head, "\shortmid"{marking}, dotted, from=4-5, to=3-7]
    \arrow[no head, from=3-7, to=4-8]
    \arrow[no head, from=4-8, to=3-9]
    \arrow["\shortmid"{marking}, curve={height=-18pt}, no head, from=1-5, to=3-9]
    \arrow["\beta"{description}, curve={height=6pt}, no head, from=6-5, to=4-8]
    \arrow["\beta"{description}, no head, from=4-5, to=6-5]
    \arrow["\beta"{description}, curve={height=6pt}, no head, from=4-3, to=6-5]
    \arrow[curve={height=6pt}, no head, from=1-5, to=3-2]
    \arrow[no head, from=1-5, to=3-4]
    \arrow[no head, from=1-5, to=3-7]
  \end{tikzcd}\]
\end{lem}
\begin{proof}
  We follow the same argument as for \reflem{lem:positive_parenthesis_completion}.
  By lozenge completion, we may find a unique $c_1$ completing the lozenge $(b,c_0,d_0,c_1)$.
  If $d_0 \prec^{-} c_1 \prec^{+} b$, then we are done taking $p=0$.
  Else, we may find a unique $d_1$ completing the lozenge $(c_1,d_0,e,d_1)$.
  Because of the sign rule, $e \prec^{\beta} d_1 \prec^{-} c_1$.
  Then we continue this process starting from $c_1$.
  It must finish because of \hyperref[prop:DFC_strictness_1]{strictness}.
\end{proof}

\noindent The two previous lemmas yield the following:

\begin{lem}\label{lem:signed_hexagon_property}
  For every hexagon as in \ref{fig:signed_hexagon},
  there is -- up to potentially exchanging $c$ and $c'$ -- a (maybe trivial) lower path from $c$ to $c'$ as in \ref{fig:signed_hexagon_filling} below:
  \vspace{-20pt}
  \begin{figure}[H]
    \begin{minipage}{170pt}
      \[\begin{tikzcd}[row sep = 15pt]
        & b \\
        c && {c'} \\
        d && {d'} \\
        & e
        \arrow[no head, from=1-2, to=2-1]
        \arrow["\alpha"', no head, from=2-1, to=3-1]
        \arrow["\beta"', no head, from=3-1, to=4-2]
        \arrow[no head, from=1-2, to=2-3]
        \arrow["{\alpha'}", no head, from=2-3, to=3-3]
        \arrow["\beta", no head, from=3-3, to=4-2]
      \end{tikzcd}\]
      \captionof{figure}{Hexagon}
      \label{fig:signed_hexagon}
    \end{minipage}
    \begin{minipage}{200pt}
      \[\begin{tikzcd}[row sep = 15pt, column sep = small]
        &&&& b \\
        \\
        c && \bullet &&& \bullet && {c'} \\
        d & \bullet && \bullet &&& \bullet & {d'} \\
        \\
        &&&& e
        \arrow["\alpha"', no head, from=3-1, to=4-1]
        \arrow[curve={height=18pt}, no head, from=1-5, to=3-1]
        \arrow[curve={height=6pt}, no head, from=1-5, to=3-3]
        \arrow[no head, from=1-5, to=3-6]
        \arrow[curve={height=-12pt}, no head, from=1-5, to=3-8]
        \arrow[no head, from=3-1, to=4-2]
        \arrow["\shortmid"{marking}, no head, from=4-2, to=3-3]
        \arrow[no head, from=3-3, to=4-4]
        \arrow[dotted, "\shortmid"{marking}, no head, from=4-4, to=3-6]
        \arrow[no head, from=3-6, to=4-7]
        \arrow["\shortmid"{marking}, no head, from=4-7, to=3-8]
        \arrow["\beta"{description}, curve={height=18pt}, no head, from=4-1, to=6-5]
        \arrow["\beta"{description}, curve={height=6pt}, no head, from=4-2, to=6-5]
        \arrow["\beta"{description}, no head, from=4-4, to=6-5]
        \arrow["\beta"{description}, no head, from=4-7, to=6-5]
        \arrow["{\alpha'}", no head, from=3-8, to=4-8]
        \arrow["\beta"{description}, curve={height=-12pt}, no head, from=4-8, to=6-5]
      \end{tikzcd}\]
      \captionof{figure}{Filled hexagon}
      \label{fig:signed_hexagon_filling}
    \end{minipage}
  \end{figure}
\end{lem}
\begin{proof}
  Suppose that we have such a hexagon, then using either \reflem{lem:positive_parenthesis_completion} or \reflem{lem:negative_parenthesis_completion}, we may produce two paths as described in the those lemmas, starting from the left side and the right side of the hexagon, respectively.
  But those two paths must finish at the same source of $\gamma(b)$. Indeed, because of uniqueness of lozenge completion, there is at most one $d''$ with $e \prec^{\beta} d'' \prec^{-} \gamma(b)$.
  Hence, both paths are ascending paths in $\delta(b)$, reaching the same leaf.
  More precisely, in the tree structure of $\delta(b)$, they both are of the form:
  $$
  c_0 \branch{d_0} c_1 \branch{d_1} \cdots \branch{d_p} c_p \branch{d_{p+1} = d''}
  $$
  Because of the tree structure, one of those paths must be an extension of the other one, which yields the result. 
\end{proof}

\begin{defin}[zig-zag]\label{defin:zig_zag}
  A \emph{zig-zag} from $c$ to $c'$ is the data of a sequence as follows:
  $$
  c = c_0 \succ^{\alpha_0} d_0 \prec^{-\alpha_0} c_1 \succ^{\alpha_1} d_1 \prec^{-\alpha_1} c_2 \succ^{\alpha_2} \cdots \succ^{\alpha_{p-1}}d_{p-1}\prec^{-\alpha_{p-1}}c_p = c'
  $$
  Such a zig-zag is said to be \emph{simple} whenever $\forall i,\, d_i \neq d_{i+1}$. \\
  It is said to be \emph{non-trivial} if $p>0$.
  If $\forall i,\, c_i \prec^{-} b$, the zig-zag will be called a $\delta(b)$\emph{-zig-zag}. \\
  Notice that because of the tree-structure on $\delta(b)$, if $c,\,c'\prec^{-}b$ there is a unique simple $\delta(b)$-zig-zag between $c$ and $c'$.
  In term of rooted trees, a simple zig-zag is a sequence of the following form:
  $$
  c_0 \branched{d_0} c_1 \branched{d_1} \cdots \branched{d_{r-1}} c_r \branch{d_r} \cdots \branch{d_{p-1}} c_{p-1} \branch{d_p} c_p
  $$
  where no two triplets are the same.
\end{defin}

\begin{prop}[Hexagon property]\label{prop:DFC_hexagon}
  For every hexagon 
  
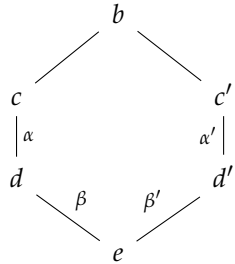
\[\begin{tikzcd}[row sep = 15pt]
    & b \\
    c && {c'} \\
    d && {d'} \\
    & e
    \arrow["\beta", no head, from=3-1, to=4-2]
    \arrow["{\beta'}"', no head, from=3-3, to=4-2]
    \arrow["\alpha", no head, from=2-1, to=3-1]
    \arrow["{\alpha'}"', no head, from=2-3, to=3-3]
    \arrow[no head, from=1-2, to=2-1]
    \arrow[no head, from=1-2, to=2-3]
  \end{tikzcd}\]
  Either $c = c'$, or (potentially by exchanging the role of $c$ and $c'$)
  there is a non-trivial simple $\delta(b)$-zig-zag as below:
  \[\begin{tikzcd}[column sep = 2pt]
    &&&&&& b \\
    \\
    {c = c_0} && {c_1} && {c_r} && {c_{r+1}} && {c_{r+2}} && {c_{p-1}} && {c_p = c'} \\
    & {d_0} && {d_1} && {d_r} && {d_{r+1}} && {d_{p-1}} && {d_{p-1}} \\
    \\
    &&&&&& e
    \arrow["\shortmid"{marking}, no head, from=3-1, to=4-2]
    \arrow[no head, from=3-3, to=4-2]
    \arrow["\shortmid"{marking}, no head, from=3-3, to=4-4]
    \arrow["\shortmid"{marking}, no head, from=3-5, to=4-6]
    \arrow[no head, from=3-7, to=4-6]
    \arrow[no head, from=3-7, to=4-8]
    \arrow["\shortmid"{marking}, no head, from=4-8, to=3-9]
    \arrow["\shortmid"{marking}, no head, from=4-10, to=3-11]
    \arrow[no head, from=3-11, to=4-12]
    \arrow["\shortmid"{marking}, no head, from=4-12, to=3-13]
    \arrow[curve={height=-18pt}, no head, from=3-1, to=1-7]
    \arrow[curve={height=-6pt}, no head, from=3-3, to=1-7]
    \arrow[no head, from=3-5, to=1-7]
    \arrow[no head, from=3-7, to=1-7]
    \arrow[no head, from=1-7, to=3-9]
    \arrow[curve={height=-6pt}, no head, from=1-7, to=3-11]
    \arrow[curve={height=-12pt}, no head, from=1-7, to=3-13]
    \arrow["\beta"{description}, no head, from=4-6, to=6-7]
    \arrow["{-\beta}"{description}, no head, from=4-8, to=6-7]
    \arrow["{-\beta}"{description}, curve={height=-6pt}, no head, from=4-10, to=6-7]
    \arrow["{-\beta}"{description}, curve={height=-12pt}, no head, from=4-12, to=6-7]
    \arrow["\beta"{description}, curve={height=6pt}, no head, from=4-4, to=6-7]
    \arrow["\beta"{description}, curve={height=12pt}, no head, from=4-2, to=6-7]
    \arrow[dotted, no head, from=4-4, to=3-5]
    \arrow[dotted, no head, from=3-9, to=4-10]
  \end{tikzcd}\]
\end{prop}
\begin{proof}
  \reflem{lem:signed_hexagon_property} shows the case $\beta = \beta'$.
  We should now consider the case where we have a hexagon as above, with $\beta' = -\beta$.
  Starting from the left part of the hexagon $e\prec^{\beta}d\prec^{\alpha}c\prec^{-}b$, we may either construct a sequence as in \ref{fig:hexagon_prop_first_case}, or as in \ref{fig:hexagon_prop_second_case} below (with $p\ge 0$):
  \begin{figure}[H]
    \begin{minipage}{.5\textwidth}
      \[\begin{tikzcd}[row sep = 15pt, column sep = tiny]
        &&& b \\
        \\
        {c = c_0} && {c_1} &&& {c_p} \\
        d & {d_0} &&& {d_{p-1}} && {d_p} \\
        \\
        &&& e
        \arrow["\alpha"', no head, from=3-1, to=4-1]
        \arrow["\shortmid"{marking}, no head, from=3-1, to=4-2]
        \arrow[no head, from=4-2, to=3-3]
        \arrow[no head, from=3-6, to=4-7]
        \arrow[curve={height=12pt}, no head, from=1-4, to=3-1]
        \arrow[no head, from=1-4, to=3-3]
        \arrow[no head, from=1-4, to=3-6]
        \arrow["\beta"{description}, curve={height=12pt}, no head, from=4-1, to=6-4]
        \arrow["{-\beta}"{description}, curve={height=-12pt}, no head, from=4-7, to=6-4]
        \arrow["\beta"{description}, no head, from=4-2, to=6-4]
        \arrow[dotted, no head, from=3-3, to=4-5]
        \arrow[no head, from=4-5, to=3-6]
        \arrow["\beta"{description}, no head, from=4-5, to=6-4]
      \end{tikzcd}\]
      \captionof{figure}{First case}
      \label{fig:hexagon_prop_first_case}
    \end{minipage}
    \begin{minipage}{.5\textwidth}
      \[\begin{tikzcd}[row sep = 15pt, column sep = tiny]
        &&& b \\
        \\
        {c = c_0} && {c_1} &&& {c_p} && {\gamma(b)} \\
        d & {d_0} &&& {d_{p-1}} && {d_p} \\
        \\
        &&& e
        \arrow["\alpha"', no head, from=3-1, to=4-1]
        \arrow["\shortmid"{marking}, no head, from=3-1, to=4-2]
        \arrow[no head, from=4-2, to=3-3]
        \arrow["\shortmid"{marking}, no head, from=3-6, to=4-7]
        \arrow[curve={height=12pt}, no head, from=1-4, to=3-1]
        \arrow[no head, from=1-4, to=3-3]
        \arrow[no head, from=1-4, to=3-6]
        \arrow["\beta"{description}, curve={height=12pt}, no head, from=4-1, to=6-4]
        \arrow["\beta"{description}, curve={height=-12pt}, no head, from=4-7, to=6-4]
        \arrow["\beta"{description}, no head, from=4-2, to=6-4]
        \arrow[dotted, no head, from=3-3, to=4-5]
        \arrow[no head, from=4-5, to=3-6]
        \arrow["\beta"{description}, no head, from=4-5, to=6-4]
        \arrow["\shortmid"{marking}, no head, from=4-7, to=3-8]
        \arrow["\shortmid"{marking}, curve={height=-18pt}, no head, from=1-4, to=3-8]
      \end{tikzcd}\]
      \captionof{figure}{Second case}
      \label{fig:hexagon_prop_second_case}
    \end{minipage}
  \end{figure}
  To see this we fill lozenges from left to right as follows:
  First, we let $d_0 := d$ if $\alpha = +$, and in the other case, we find $d_0$ as the unique one completing the lozenge $(c_0,\, d,\, e,\, d_0)$.
  If $e\prec^{-\beta} d_0\prec^{-} c_0$, then we are in the first case, with $p = 0$.
  Else we find the unique $c_1$ filling the lozenge $(b,\, c_0,\, d_0,\, c_1)$.
  If $d_0 \prec^{+}c_1\prec^{+}b$, then we are in the second case with $p = 0$.
  Else, we continue the same process, filling the lozenge $(c_1,\,d_0,\,e,\,d_1)$ \textit{etc.}
  This process must finish because of \hyperref[defin:PFS]{strictness}.

  Now, if we end up in the first case, we may use \reflem{lem:signed_hexagon_property} to find a path between $c_p$ and $c'$, and conclude.
  If we end up in the second case, we repeat the same argument with the right parenthesis $e \prec^{-\beta} d' \prec^{\alpha'} c' \prec^{-} b$.
  We thus find a path, either with the shape of \ref{fig:hexagon_prop_first_case_bis}, or with the shape of \ref{fig:hexagon_prop_second_case_bis} below (with $q\ge0$):
  \begin{figure}[H]
    \begin{minipage}{.5\textwidth}
      \[\begin{tikzcd}[row sep = 15pt, column sep = tiny]
        &&& b \\
        \\
        {c' = c_0'} && {c_1'} &&& {c_q'} \\
        d' & {d_0'} &&& {d_{q-1}'} && {d_q'} \\
        \\
        &&& e
        \arrow["\alpha'"', no head, from=3-1, to=4-1]
        \arrow["\shortmid"{marking}, no head, from=3-1, to=4-2]
        \arrow[no head, from=4-2, to=3-3]
        \arrow[no head, from=3-6, to=4-7]
        \arrow[curve={height=12pt}, no head, from=1-4, to=3-1]
        \arrow[no head, from=1-4, to=3-3]
        \arrow[no head, from=1-4, to=3-6]
        \arrow["-\beta"{description}, curve={height=12pt}, no head, from=4-1, to=6-4]
        \arrow["\beta"{description}, curve={height=-12pt}, no head, from=4-7, to=6-4]
        \arrow["-\beta"{description}, no head, from=4-2, to=6-4]
        \arrow[dotted, no head, from=3-3, to=4-5]
        \arrow[no head, from=4-5, to=3-6]
        \arrow["-\beta"{description}, no head, from=4-5, to=6-4]
      \end{tikzcd}\]
      \captionof{figure}{First case}
      \label{fig:hexagon_prop_first_case_bis}
    \end{minipage}
    \begin{minipage}{.5\textwidth}
      \[\begin{tikzcd}[row sep = 15pt, column sep = tiny]
        &&& b \\
        \\
        {c' = c_0'} && {c_1'} &&& {c_q'} && {\gamma(b)} \\
        d' & {d_0'} &&& {d_{q-1}'} && {d_q'} \\
        \\
        &&& e
        \arrow["\alpha'"', no head, from=3-1, to=4-1]
        \arrow["\shortmid"{marking}, no head, from=3-1, to=4-2]
        \arrow[no head, from=4-2, to=3-3]
        \arrow["\shortmid"{marking}, no head, from=3-6, to=4-7]
        \arrow[curve={height=12pt}, no head, from=1-4, to=3-1]
        \arrow[no head, from=1-4, to=3-3]
        \arrow[no head, from=1-4, to=3-6]
        \arrow["-\beta"{description}, curve={height=12pt}, no head, from=4-1, to=6-4]
        \arrow["-\beta"{description}, curve={height=-12pt}, no head, from=4-7, to=6-4]
        \arrow["-\beta"{description}, no head, from=4-2, to=6-4]
        \arrow[dotted, no head, from=3-3, to=4-5]
        \arrow[no head, from=4-5, to=3-6]
        \arrow["-\beta"{description}, no head, from=4-5, to=6-4]
        \arrow["\shortmid"{marking}, no head, from=4-7, to=3-8]
        \arrow["\shortmid"{marking}, curve={height=-18pt}, no head, from=1-4, to=3-8]
      \end{tikzcd}\]
      \captionof{figure}{Second case}
      \label{fig:hexagon_prop_second_case_bis}
    \end{minipage}
  \end{figure}
  In the first case, we may conclude using \reflem{lem:signed_hexagon_property}.
  In the second case, we end up with the following configuration:
  
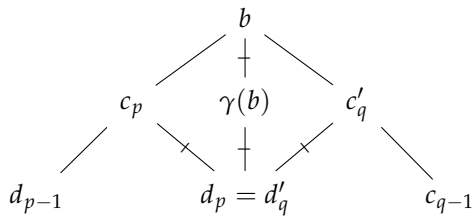
\[\begin{tikzcd}[row sep = 15pt, column sep = small]
    && b \\
    & {c_p} & {\gamma(b)} & {c_q'} \\
    {d_{p-1}} && {d_p = d_q'} && {c_{q-1}}
    \arrow[no head, from=1-3, to=2-2]
    \arrow["\shortmid"{marking}, no head, from=1-3, to=2-3]
    \arrow["\shortmid"{marking}, no head, from=2-3, to=3-3]
    \arrow["\shortmid"{marking}, no head, from=2-2, to=3-3]
    \arrow[no head, from=1-3, to=2-4]
    \arrow["\shortmid"{marking}, no head, from=2-4, to=3-3]
    \arrow[no head, from=2-4, to=3-5]
    \arrow[no head, from=3-1, to=2-2]
  \end{tikzcd}\]
  So, by uniqueness of lozenge completion, $c_p = c_q'$ and we have the following zig-zag between $c$ and $c'$:
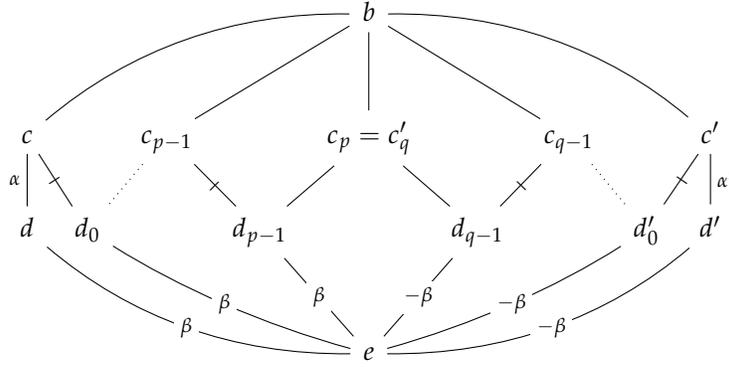
  \[\begin{tikzcd}[row sep = 15pt, column sep = tiny]
    &&&& b \\
    \\
    c && {c_{p-1}} && {c_p=c_q'} && {c_{q-1}} && {c'} \\
    d & {d_0} && {d_{p-1}} && {d_{q-1}} && {d_0'} & {d'} \\
    \\
    &&&& e
    \arrow[no head, from=3-5, to=4-4]
    \arrow[no head, from=3-5, to=4-6]
    \arrow["\shortmid"{marking}, no head, from=3-7, to=4-6]
    \arrow["\shortmid"{marking}, no head, from=3-3, to=4-4]
    \arrow[no head, from=1-5, to=3-3]
    \arrow[no head, from=1-5, to=3-5]
    \arrow[no head, from=1-5, to=3-7]
    \arrow[dotted, no head, from=3-3, to=4-2]
    \arrow[dotted, no head, from=3-7, to=4-8]
    \arrow["\shortmid"{marking}, no head, from=4-8, to=3-9]
    \arrow["\shortmid"{marking}, no head, from=3-1, to=4-2]
    \arrow["\alpha"', no head, from=3-1, to=4-1]
    \arrow["{\alpha'}", no head, from=3-9, to=4-9]
    \arrow[curve={height=18pt}, no head, from=1-5, to=3-1]
    \arrow[curve={height=-18pt}, no head, from=1-5, to=3-9]
    \arrow["\beta"{description}, curve={height=18pt}, no head, from=4-1, to=6-5]
    \arrow["{-\beta}"{description}, curve={height=-18pt}, no head, from=4-9, to=6-5]
    \arrow["\beta"{description}, no head, from=4-4, to=6-5]
    \arrow["{-\beta}"{description}, no head, from=4-6, to=6-5]
    \arrow["{-\beta}"{description}, curve={height=-6pt}, no head, from=4-8, to=6-5]
    \arrow["\beta"{description}, curve={height=6pt}, no head, from=4-2, to=6-5]
  \end{tikzcd}\]
\end{proof}

\begin{prop}[pencil linearity]\label{prop:DFC_pencil_linearity}
  $C$ satisfies the axiom of \hyperref[defin:PFS]{pencil linearity}.
  That is:\\ $\forall k>0,\,\forall e\in C_{k-1},\,\forall \beta \in\{+,\,-\},\quad\left\{d\in C_k \mid e\prec^{\beta}d\right\}$ is linearly ordered by $<^{+}$.
\end{prop}
\begin{proof}
  Suppose $d\succ^{\beta}e\prec^{\beta}d'$, with $d \neq d'$.
  We know that every element is always the source or the target of some element in $\Lambda$.
  Because elements in $\Lambda$ are sources of an iterated target of $\omega$, we may find a hexagone as follows:
  
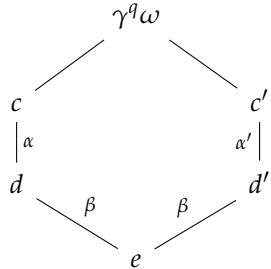
\[\begin{tikzcd}[row sep = 15pt]
    & \gamma^{q}\omega \\
    c && {c'} \\
    d && {d'} \\
    & e
    \arrow["\beta", no head, from=3-1, to=4-2]
    \arrow["{\beta}"', no head, from=3-3, to=4-2]
    \arrow["\alpha", no head, from=2-1, to=3-1]
    \arrow["{\alpha'}"', no head, from=2-3, to=3-3]
    \arrow[no head, from=1-2, to=2-1]
    \arrow[no head, from=1-2, to=2-3]
  \end{tikzcd}\]
  If $c = c'$, then by the sign rule $\alpha = - \alpha'$, hence $d$ and $d'$ are $<^{+}$-comparable and we are done.
  In what follows, we suppose $c \neq c'$.
  So -- up to potentially exchanging $c$ and $c'$ --
  there is a simple non-trivial zig-zag with $d_0 = \gamma(c_0)$ as follows:
  
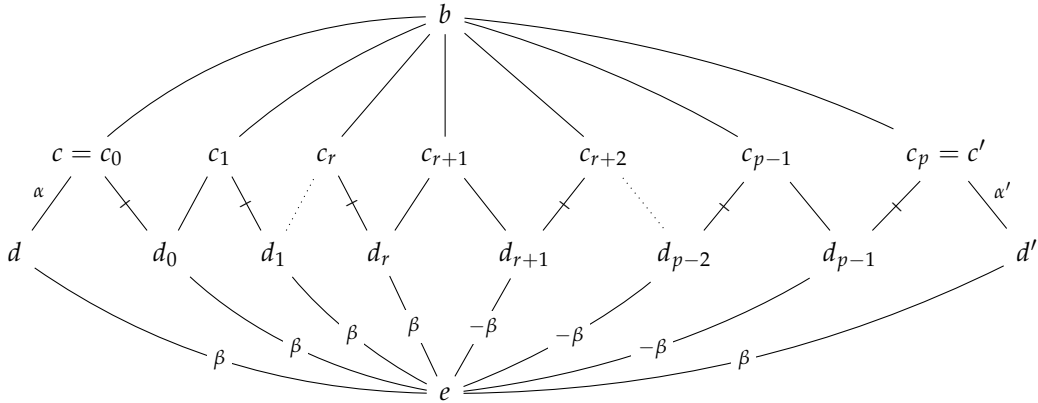
\[\begin{tikzcd}[column sep = 2pt]
    &&&&&&& b \\
    \\
    & {c = c_0} && {c_1} && {c_r} && {c_{r+1}} && {c_{r+2}} && {c_{p-1}} && {c_p = c'} \\
    d && {d_0} && {d_1} && {d_r} && {d_{r+1}} && {d_{p-2}} && {d_{p-1}} && {d'} \\
    \\
    &&&&&&& e
    \arrow["\shortmid"{marking}, no head, from=3-2, to=4-3]
    \arrow[no head, from=3-4, to=4-3]
    \arrow["\shortmid"{marking}, no head, from=3-4, to=4-5]
    \arrow["\shortmid"{marking}, no head, from=3-6, to=4-7]
    \arrow[no head, from=3-8, to=4-7]
    \arrow[no head, from=3-8, to=4-9]
    \arrow["\shortmid"{marking}, no head, from=4-9, to=3-10]
    \arrow["\shortmid"{marking}, no head, from=4-11, to=3-12]
    \arrow[no head, from=3-12, to=4-13]
    \arrow["\shortmid"{marking}, no head, from=4-13, to=3-14]
    \arrow[curve={height=-18pt}, no head, from=3-2, to=1-8]
    \arrow[curve={height=-6pt}, no head, from=3-4, to=1-8]
    \arrow[no head, from=3-6, to=1-8]
    \arrow[no head, from=3-8, to=1-8]
    \arrow[no head, from=1-8, to=3-10]
    \arrow[curve={height=-6pt}, no head, from=1-8, to=3-12]
    \arrow[curve={height=-12pt}, no head, from=1-8, to=3-14]
    \arrow["\beta"{description}, no head, from=4-7, to=6-8]
    \arrow["{-\beta}"{description}, no head, from=4-9, to=6-8]
    \arrow["{-\beta}"{description}, curve={height=-6pt}, no head, from=4-11, to=6-8]
    \arrow["{-\beta}"{description}, curve={height=-12pt}, no head, from=4-13, to=6-8]
    \arrow["\beta"{description}, curve={height=6pt}, no head, from=4-5, to=6-8]
    \arrow["\beta"{description}, curve={height=12pt}, no head, from=4-3, to=6-8]
    \arrow["{\alpha'}", no head, from=3-14, to=4-15]
    \arrow["\alpha"', no head, from=3-2, to=4-1]
    \arrow["\beta"{description}, curve={height=18pt}, no head, from=4-1, to=6-8]
    \arrow["\beta"{description}, curve={height=-18pt}, no head, from=4-15, to=6-8]
    \arrow[dotted, no head, from=4-5, to=3-6]
    \arrow[dotted, no head, from=3-10, to=4-11]
  \end{tikzcd}\]
  In this situation, distinguishing on the sign of $\alpha'$, we see on the rightmost lozenge that necessarily $r = p-1$.
  Hence the zig-zag is a path:
  
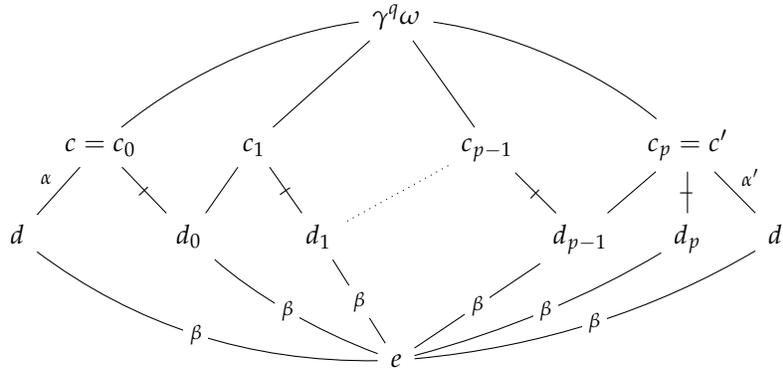
\[\begin{tikzcd}[column sep = tiny, row sep = 15pt]
    &&&&& {\gamma^q\omega} \\
    \\
    & {c = c_0} && {c_1} &&& {c_{p-1}} && {c_p = c'} \\
    d && {d_0} && {d_1} &&& {d_{p-1}} & {d_p} & {d'} \\
    \\
    &&&&& e
    \arrow[dotted, no head, from=3-7, to=4-5]
    \arrow[no head, from=3-4, to=4-3]
    \arrow["\shortmid"{marking}, no head, from=3-2, to=4-3]
    \arrow["\shortmid"{marking}, no head, from=3-4, to=4-5]
    \arrow["\shortmid"{marking}, no head, from=3-7, to=4-8]
    \arrow[no head, from=4-8, to=3-9]
    \arrow["\alpha"', no head, from=3-2, to=4-1]
    \arrow["{\alpha'}", no head, from=3-9, to=4-10]
    \arrow["\beta"{description}, curve={height=18pt}, no head, from=4-1, to=6-6]
    \arrow["\beta"{description}, curve={height=6pt}, no head, from=4-3, to=6-6]
    \arrow["\beta"{description}, no head, from=4-5, to=6-6]
    \arrow["\beta"{description}, no head, from=4-8, to=6-6]
    \arrow["\beta"{description}, curve={height=-12pt}, no head, from=4-10, to=6-6]
    \arrow[curve={height=12pt}, no head, from=1-6, to=3-2]
    \arrow[no head, from=1-6, to=3-4]
    \arrow[no head, from=1-6, to=3-7]
    \arrow[curve={height=-12pt}, no head, from=1-6, to=3-9]
    \arrow["\shortmid"{marking}, no head, from=3-9, to=4-9]
    \arrow["\beta"{description}, curve={height=-6pt}, no head, from=4-9, to=6-6]
  \end{tikzcd}\]
  Either $\alpha = +$ and $d = d_0$, or $\alpha= -$ and $d\trigup{}d_0$.
  Also, either $\alpha' = +$ and $d = d_p$, or $\alpha' = -$, and because of the uniqueness of lozenge completion, $d' = d_{p-1}$.
  Summing up all cases, we always have an upper path from $c$ to $c'$ (in the tree structure of $\gamma^q\omega$, and of the shape depicted above).
  Hence it proves that all faces having a given $e$ as a source are linearly ordered for $<^{+}$, and that all faces having $e$ as a target also are linearly ordered by $<^{+}$.
\end{proof}

\para We shall finally focus on the remaining axiom of positive opetopes: \hyperref[defin:PFS]{disjointness}.
We will need some lemmas, and we will split the proof into two parts.

\begin{lem}\label{lem:DFC_path_along_the_border_of_a_cell}
  If $d\prec^{-}c$ with $d\in C_{\ge 1}$, then there is a (unique) path as follows.
  
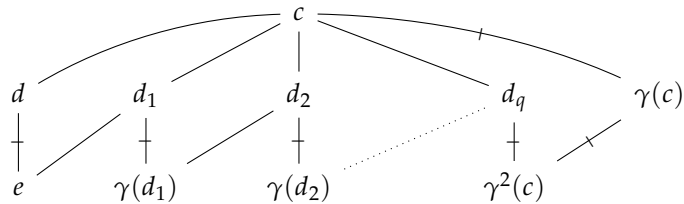
\[\begin{tikzcd}[row sep = 15pt]
    && c \\
    d & {d_1} & {d_2} && {d_q} & {\gamma(c)} \\
    e & {\gamma(d_1)} & {\gamma(d_2)} && {\gamma^2(c)}
    \arrow[no head, from=3-1, to=2-2]
    \arrow["\shortmid"{marking}, no head, from=2-1, to=3-1]
    \arrow["\shortmid"{marking}, no head, from=2-2, to=3-2]
    \arrow[no head, from=2-3, to=3-2]
    \arrow["\shortmid"{marking}, no head, from=2-3, to=3-3]
    \arrow[no head, dotted, from=3-3, to=2-5]
    \arrow["\shortmid"{marking}, no head, from=2-5, to=3-5]
    \arrow["\shortmid"{marking}, no head, from=2-6, to=3-5]
    \arrow[curve={height=10pt}, no head, from=1-3, to=2-1]
    \arrow[no head, from=2-2, to=1-3]
    \arrow[no head, from=1-3, to=2-3]
    \arrow[no head, from=1-3, to=2-5]
    \arrow["\shortmid"{marking}, curve={height=-10pt}, no head, from=1-3, to=2-6]
  \end{tikzcd}\]
  In fact this is the unique path from $d$ to $\rho(\delta(c))$ in the tree structure on $\delta(c)$.
\end{lem}
\begin{proof}
  This is seen by completing lozenges from left to right until
  coming accross one with the shape
  
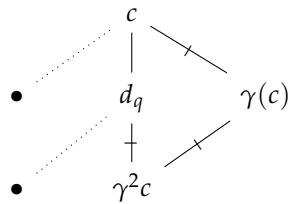
\[\begin{tikzcd}[row sep = 15pt]
    & c \\
    \bullet & {d_q} & {\gamma(c)} \\
    \bullet & {\gamma^2c}
    \arrow[dotted, no head, from=3-1, to=2-2]
    \arrow[no head, from=1-2, to=2-2]
    \arrow[dotted, no head, from=2-1, to=1-2]
    \arrow["\shortmid"{marking}, no head, from=2-2, to=3-2]
    \arrow["\shortmid"{marking}, no head, from=1-2, to=2-3]
    \arrow["\shortmid"{marking}, no head, from=2-3, to=3-2]
  \end{tikzcd}\]
  which ends the path.
  The process must finish because there is no infinite branch in the tree structure of $\delta(c)$.
  Moreover the path obtained this way is unique because all lozenge completions are.
\end{proof}

This lemma yields the following one, by concatenating paths.
\begin{lem}\label{lem:DFC_path_along_the_border_of_several_cells}
  If there is an upper path: $d = d_0\prec^{-}c_1\succ^{+}\gamma(c_1)\prec^{-}\cdots\prec^{-}c_q\succ^{+}\gamma(c_q) = d'$,
  then there is a path:
  
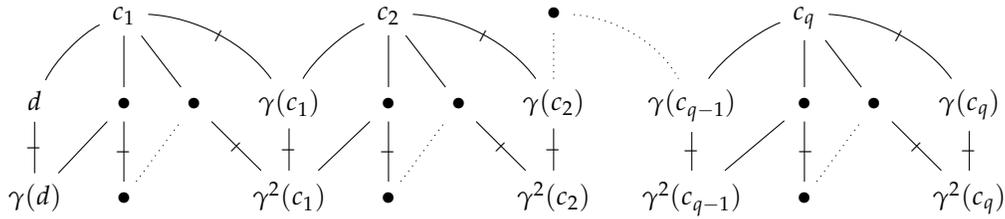
\[\begin{tikzcd}[column sep = 10pt, row sep = 15pt]
    & {c_1} &&& {c_2} && \bullet && {c_q} \\
    d & \bullet & \bullet & {\gamma(c_1)} & \bullet & \bullet & {\gamma(c_2)} & {\gamma(c_{q-1})} & \bullet & \bullet & {\gamma(c_q)} \\
    {\gamma(d)} & \bullet && {\gamma^2(c_1)} & \bullet && {\gamma^2(c_2)} & {\gamma^2(c_{q-1})} & \bullet && {\gamma^2(c_q)}
    \arrow[no head, from=3-1, to=2-2]
    \arrow["\shortmid"{marking}, no head, from=2-1, to=3-1]
    \arrow["\shortmid"{marking}, no head, from=2-2, to=3-2]
    \arrow[curve={height=6pt}, no head, from=1-2, to=2-1]
    \arrow[no head, from=2-2, to=1-2]
    \arrow[dotted, no head, from=3-2, to=2-3]
    \arrow[no head, from=1-2, to=2-3]
    \arrow["\shortmid"{marking}, curve={height=-12pt}, no head, from=1-2, to=2-4]
    \arrow["\shortmid"{marking}, no head, from=2-3, to=3-4]
    \arrow["\shortmid"{marking}, no head, from=2-4, to=3-4]
    \arrow[curve={height=-6pt}, no head, from=2-4, to=1-5]
    \arrow[no head, from=1-5, to=2-5]
    \arrow[no head, from=3-4, to=2-5]
    \arrow["\shortmid"{marking}, no head, from=2-5, to=3-5]
    \arrow[no head, from=1-5, to=2-6]
    \arrow[dotted, no head, from=3-5, to=2-6]
    \arrow["\shortmid"{marking}, no head, from=2-6, to=3-7]
    \arrow["\shortmid"{marking}, curve={height=-12pt}, no head, from=1-5, to=2-7]
    \arrow["\shortmid"{marking}, no head, from=2-7, to=3-7]
    \arrow[dotted, no head, from=2-7, to=1-7]
    \arrow[curve={height=-12pt}, dotted, no head, from=1-7, to=2-8]
    \arrow["\shortmid"{marking}, no head, from=2-8, to=3-8]
    \arrow[curve={height=-6pt}, no head, from=2-8, to=1-9]
    \arrow[no head, from=2-9, to=1-9]
    \arrow[no head, from=2-10, to=1-9]
    \arrow["\shortmid"{marking}, curve={height=12pt}, no head, from=2-11, to=1-9]
    \arrow["\shortmid"{marking}, no head, from=2-11, to=3-11]
    \arrow["\shortmid"{marking}, no head, from=2-10, to=3-11]
    \arrow["\shortmid"{marking}, no head, from=3-9, to=2-9]
    \arrow[no head, from=3-8, to=2-9]
    \arrow[dotted, no head, from=3-9, to=2-10]
  \end{tikzcd}\]
  hence yielding an upper path from $\gamma(d)$ to $\gamma(d')$
\end{lem}
Now using \reflem{lem:DFC_path_along_the_border_of_several_cells}, we are able to prove the first half of \hyperref[defin:PFS]{disjointness}:

\begin{prop}\label{prop:DFC_disjointness_1}
  If $k>0$, and $d,\,d'\in C_k$ then it cannot be the case that $d<^{+}d'$ and $d'<^{-}d$.
\end{prop}
\begin{proof}
  If $d<^{+}d'$ and $d'<^{-}d$, then we have $\gamma(d)\le^{+}\gamma(d')$ by \reflem{lem:DFC_path_along_the_border_of_several_cells}, and $\gamma(d')<^{+}\gamma(d)$ by the second hypothesis (because $d'<^{-}d \Rightarrow \gamma(d') <^{+} \gamma(d)$).
  Thus, we have a cycle $\gamma(d) <^{+} \gamma(d)$, which is absurd by \hyperref[prop:DFC_strictness_1]{strictness}.
\end{proof}

\para We will now focus on the second half by handling the case where $d<^{+}d'$, and $d<^{-}d'$.
We begin with another lemma.

\begin{lem}\label{lem:DFC_no_upper_path_in_cell}
  There is no (non-trivial) upper path between two sources of a common cell.
\end{lem}
\begin{proof}
  First, if the common cell is $\omega$, one of the sources should also be the target of $\omega$, which is absurd.
  If the common cell $a_0$ is of codimension $1$, assuming the upper path $b_p \prec^{-} a_p \triangleleft^{-} a_{p-1} \triangleleft^{-} \cdots \triangleleft^{-} a_1 \succ^{+} b_0$ takes place in $\Lambda$, we have a diagram as follows:
  
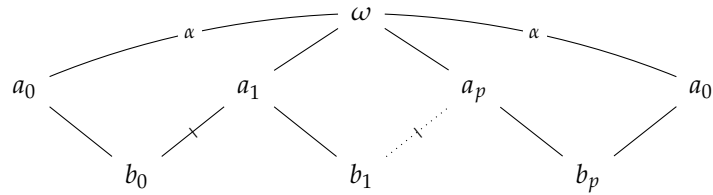
\[\begin{tikzcd}[row sep = 15pt]
    &&& \omega \\
    {a_0} && {a_1} && {a_p} && {a_0} \\
    & {b_0} && {b_1} && {b_p}
    \arrow[no head, from=3-6, to=2-7]
    \arrow[no head, from=3-6, to=2-5]
    \arrow[no head, from=3-4, to=2-3]
    \arrow[no head, from=3-2, to=2-1]
    \arrow["\shortmid"{marking}, no head, from=2-3, to=3-2]
    \arrow["\shortmid"{marking}, dotted, no head, from=3-4, to=2-5]
    \arrow["\alpha"{description}, curve={height=8pt}, no head, from=1-4, to=2-1]
    \arrow[no head, from=1-4, to=2-3]
    \arrow[no head, from=1-4, to=2-5]
    \arrow["\alpha"{description}, curve={height=-8pt}, no head, from=1-4, to=2-7]
  \end{tikzcd}\]
  and we cannot choose $\alpha$ consistently.

  From now on, we will assume that the dimension of the common cell is $n-k$ with $k\ge2$.

  Notice that we have the following corollary of the \hyperref[prop:DFC_hexagon]{hexagon property}:
  If there is a hexagon as in \ref{fig:no_upper_path_hexagon} below,
  \begin{figure}[H]
    \begin{minipage}{.5\textwidth}
      \[\begin{tikzcd}[row sep = 15pt, column sep = small]
        & {\gamma^q(\omega)} \\
        c && {c'} \\
        d && {d'} \\
        & e
        \arrow[no head, from=3-1, to=4-2]
        \arrow["\shortmid"{marking}, no head, from=3-3, to=4-2]
        \arrow["\alpha", no head, from=2-1, to=3-1]
        \arrow["{\alpha'}"', no head, from=2-3, to=3-3]
        \arrow[no head, from=1-2, to=2-1]
        \arrow[no head, from=1-2, to=2-3]
      \end{tikzcd}\]
      \captionof{figure}{Hexagon}
      \label{fig:no_upper_path_hexagon}
    \end{minipage}
    \begin{minipage}{.5\textwidth}
      \[\begin{tikzcd}[row sep = 15pt]
        && {\gamma^q(\omega)} \\
        {c_0} & {c_1} & {c_2} && {c_p} \\
        {d_0} & {d_1} & {d_2} && {d_p} \\
        {e_0} & {e_1} & {e_2}
        \arrow["{\alpha_0}"', no head, from=2-1, to=3-1]
        \arrow[no head, curve={height=8pt}, from=1-3, to=2-1]
        \arrow[no head, curve={height=-8pt}, from=1-3, to=2-5]
        \arrow["{\alpha_p}"', no head, from=2-5, to=3-5]
        \arrow[no head, from=3-1, to=4-1]
        \arrow["\shortmid"{marking}, no head, from=4-1, to=3-2]
        \arrow[no head, from=3-2, to=4-2]
        \arrow["\shortmid"{marking}, no head, from=4-2, to=3-3]
        \arrow[no head, from=3-3, to=4-3]
        \arrow["\shortmid"{marking}, dotted, no head, from=4-3, to=3-5]
        \arrow["{\alpha_1}"', no head, from=2-2, to=3-2]
        \arrow["{\alpha_2}"', no head, from=2-3, to=3-3]
        \arrow[no head, from=2-2, to=1-3]
        \arrow[no head, from=2-3, to=1-3]
      \end{tikzcd}\]
      \captionof{figure}{Succesive hexagons}
      \label{fig:no_upper_path_shape}
    \end{minipage}
  \end{figure}
  \noindent then there is a simple zig-zag from $c$ to $c'$ in the tree structure of $\gamma^q\omega$, where every lower element has $e$ as a codimension-$1$ subface.
  We can generalize this assertion a bit:
  When there is a shape as in \ref{fig:no_upper_path_shape} above,
  there is a zig-zag from $c_0$ to $c_p$ in the tree-structure of $\gamma^q\omega$ where lower elements successively have $e_1$, then $e_2$ ... then $e_n$ as a codimension-$1$ subface.
  It is seen by considering such zig-zags between each $c_i$ and $c_{i+1}$.
  More precisely, assuming that for all $i,\, \alpha_i = -$, the -- \textit{a priori} non simple -- zig-zag has the shape depicted in \ref{fig:huge_path}.
  \begin{sidewaysfigure}[p]
    \scalebox{0.9}{
      \begin{tikzcd}[ampersand replacement=\&,cramped,column sep=tiny, row sep = 40pt]
        \&\&\& {\gamma^q(\omega)} \&\&\&\&\&\& {\gamma^q(\omega)} \&\&\&\&\&\&\&\& {\gamma^q(\omega)} \\
        \\
        {c_0 = c_{0,0}} \& {c_{0,r_0}} \&\& {c_{0,r_0+1}} \&\& {c_{0,\,r_0+2}} \& {c_1 = c_{1,0}} \& {c_{1,r_1}} \&\& {c_{1,r_1+1}} \&\& {c_{1,r_1+2}} \& {c_2 = c_{2,0}} \&\& {c_{p-1}} \& {c_{p,r_p}} \&\& {c_{p,r_p+1}} \&\& {c_{p,r_p+2}} \& {c_p} \\
        {d_0} \& {d_{0,0}} \& {d_{0,r_0}} \&\& {d_{0,r+1}} \& {d_{0,s_0-1}} \& {d_1} \& {d_{1,0}} \& {d_{1,r_1}} \&\& {d_{1,r_1+1}} \& {d_{1,s_1-1}} \& {d_2} \&\& {d_{p-1}} \& {d_{p,0}} \& {d_{p,r_p}} \&\& {d_{p,r_p+1}} \& {d_{p,s_p-1}} \& {d_p} \\
        \\
        \&\&\& {e_0} \&\&\&\&\&\& {e_1} \&\&\&\&\&\&\&\& {e_p}
        \arrow["\shortmid"{marking}, no head, from=3-1, to=4-2]
        \arrow["\shortmid"{marking}, no head, from=3-2, to=4-3]
        \arrow[no head, from=3-4, to=4-3]
        \arrow[no head, from=3-4, to=4-5]
        \arrow["\shortmid"{marking}, no head, from=4-5, to=3-6]
        \arrow["\shortmid"{marking}, no head, from=4-6, to=3-7]
        \arrow[curve={height=-18pt}, no head, from=3-1, to=1-4]
        \arrow[no head, from=3-2, to=1-4]
        \arrow[no head, from=3-4, to=1-4]
        \arrow[no head, from=1-4, to=3-6]
        \arrow[curve={height=-18pt}, no head, from=1-4, to=3-7]
        \arrow[no head, from=4-3, to=6-4]
        \arrow["\shortmid"{marking}, no head, from=4-5, to=6-4]
        \arrow["\shortmid"{marking}, curve={height=-6pt}, no head, from=4-6, to=6-4]
        \arrow[curve={height=6pt}, no head, from=4-2, to=6-4]
        \arrow["\shortmid"{marking}, curve={height=-18pt}, no head, from=4-7, to=6-4]
        \arrow[no head, from=3-7, to=4-7]
        \arrow[no head, from=3-1, to=4-1]
        \arrow[curve={height=18pt}, no head, from=4-1, to=6-4]
        \arrow["\shortmid"{marking}, no head, from=3-7, to=4-8]
        \arrow["\shortmid"{marking}, no head, from=3-8, to=4-9]
        \arrow[no head, from=4-9, to=3-10]
        \arrow[no head, from=3-10, to=4-11]
        \arrow["\shortmid"{marking}, no head, from=4-11, to=3-12]
        \arrow["\shortmid"{marking}, no head, from=4-12, to=3-13]
        \arrow[curve={height=18pt}, no head, from=1-10, to=3-7]
        \arrow[no head, from=3-8, to=1-10]
        \arrow[no head, from=3-10, to=1-10]
        \arrow[no head, from=3-12, to=1-10]
        \arrow[curve={height=-18pt}, no head, from=1-10, to=3-13]
        \arrow[curve={height=18pt}, no head, from=4-7, to=6-10]
        \arrow["\shortmid"{marking}, curve={height=-6pt}, no head, from=4-12, to=6-10]
        \arrow[no head, from=3-13, to=4-13]
        \arrow["\shortmid"{marking}, curve={height=-18pt}, no head, from=4-13, to=6-10]
        \arrow[curve={height=6pt}, no head, from=4-8, to=6-10]
        \arrow["\shortmid"{marking}, no head, from=4-11, to=6-10]
        \arrow[no head, from=4-9, to=6-10]
        \arrow["\shortmid"{marking}, no head, from=3-15, to=4-16]
        \arrow[no head, from=3-15, to=4-15]
        \arrow["\shortmid"{marking}, no head, from=3-16, to=4-17]
        \arrow[no head, from=3-18, to=4-17]
        \arrow[no head, from=3-18, to=4-19]
        \arrow["\shortmid"{marking}, no head, from=3-20, to=4-19]
        \arrow["\shortmid"{marking}, no head, from=4-20, to=3-21]
        \arrow[no head, from=3-21, to=4-21]
        \arrow[curve={height=18pt}, no head, from=1-18, to=3-15]
        \arrow[no head, from=1-18, to=3-16]
        \arrow[no head, from=1-18, to=3-18]
        \arrow[no head, from=1-18, to=3-20]
        \arrow[curve={height=-18pt}, no head, from=1-18, to=3-21]
        \arrow[curve={height=-18pt}, no head, from=6-18, to=4-15]
        \arrow[curve={height=-6pt}, no head, from=6-18, to=4-16]
        \arrow[no head, from=6-18, to=4-17]
        \arrow["\shortmid"{marking}, no head, from=6-18, to=4-19]
        \arrow["\shortmid"{marking}, curve={height=6pt}, no head, from=6-18, to=4-20]
        \arrow["\shortmid"{marking}, curve={height=18pt}, no head, from=6-18, to=4-21]
        \arrow[curve={height=-30pt}, dotted, no head, from=3-13, to=3-15]
        \arrow[curve={height=30pt}, dotted, no head, from=4-13, to=4-15]
        \arrow[dotted, no head, from=3-20, to=4-20]
        \arrow[dotted, no head, from=3-16, to=4-16]
        \arrow[dotted, no head, from=3-12, to=4-12]
        \arrow[dotted, no head, from=3-8, to=4-8]
        \arrow[dotted, no head, from=4-2, to=3-2]
        \arrow[dotted, no head, from=4-6, to=3-6]
      \end{tikzcd}
    }
    \vspace{50pt}
    \captionof{figure}{Long zig-zag}
    \label{fig:huge_path}
  \end{sidewaysfigure}
  Then, by shortening this zig-zag whenever possible, one obtains a simple zig-zag between $c_0$ and $c_p$ made of consecutive triplets $$c_{\varphi(i),\,\sigma(i)}\succ^{\alpha_i}d_{\varphi(i),\,\sigma(i)}\prec^{-\alpha_i}c_{\varphi(i),\,\sigma(i)+1}$$ for $0\le i<l$ with $\forall i\le l-2,\,\, c_{\varphi(i),\,\sigma(i)+1} = c_{\varphi(i+1),\,\sigma(i+1)}$.\\
  Moreover, $\varphi$ is non-decreasing, $(\alpha_i)_i$ is non-increasing, and if $\varphi(i) = \varphi(i+1)$ then $\sigma(i) < \sigma(i+1)$.
  Notice that for every $i$, $e_{\varphi(i)}\prec^{-\alpha_i}d_{\varphi(i),\,\sigma(i)}$.
  
  Now suppose that we have a path as in \ref{fig:dfc_to_zpo_hasse_13} below, such that $e_p$ and $e_0$ are both sources of a same element: $e_0,\,e_p \prec^{-} d_0 \in C_{n-k}$.
  We may moreover suppose that all $d_i$'s for $i>0$ are in $\Lambda$.
  \begin{figure}[H]
    \begin{minipage}{.5\textwidth}
      \[\begin{tikzcd}[row sep = 15pt, column sep = tiny]
        {d_0} && {d_1} &&& {d_p} && {d_0} \\
        & {e_0} && {e_1} &&& {e_p}
        \arrow[no head, from=2-7, to=1-8]
        \arrow[no head, from=2-7, to=1-6]
        \arrow[no head, from=2-4, to=1-3]
        \arrow[no head, from=2-2, to=1-1]
        \arrow["\shortmid"{marking}, no head, from=1-3, to=2-2]
        \arrow["\shortmid"{marking}, dotted, no head, from=2-4, to=1-6]
      \end{tikzcd}\]
      \captionof{figure}{Hypothetical path}
      \label{fig:dfc_to_zpo_hasse_13}
    \end{minipage}
    \begin{minipage}{.5\textwidth}
      \[\begin{tikzcd}[row sep = 15pt, column sep = small]
        & {\gamma^{k-1}(\omega)} \\
        {d_0} && {d_1} \\
        & {e_0}
        \arrow[no head, from=3-2, to=2-1]
        \arrow["\shortmid"{marking}, no head, from=3-2, to=2-3]
        \arrow["\shortmid"{marking}, no head, from=2-1, to=1-2]
        \arrow[no head, from=2-3, to=1-2]
      \end{tikzcd}\]
      \captionof{figure}{Impossible lozenge}
      \label{fig:dfc_to_zpo_hasse_14}
    \end{minipage}
  \end{figure}
  If $d_0 = \gamma^k \omega$, then there is a lozenge as in \ref{fig:dfc_to_zpo_hasse_14} above,
  which is not possible by the sign rule.
  Hence $d_0$ is not $\gamma^k(\omega$), and is (by \reflem{lem:DFC_not_a_source_is_iterated_target} and \reflem{lem:DFC_source_is_source_of_Lambda}) the source of some $c_0\in\Lambda$.

  So we may extend this path with hexagons as in \ref{fig:dfc_to_zpo_hasse_15} below.
  \begin{figure}[H]
    \begin{minipage}{.5\textwidth}
      \[\begin{tikzcd}[row sep = small, column sep = tiny]
        &&& {\gamma^{k-2}\omega} \\
        \\
        {c_0} && {c_1} &&& {c_p} && {c_0} \\
        {d_0} && {d_1} &&& {d_p} && {d_0} \\
        & {e_0} && {e_1} &&& {e_p}
        \arrow[no head, from=5-7, to=4-8]
        \arrow[no head, from=5-7, to=4-6]
        \arrow[no head, from=5-4, to=4-3]
        \arrow[no head, from=5-2, to=4-1]
        \arrow["\shortmid"{marking}, no head, from=4-3, to=5-2]
        \arrow["\shortmid"{marking}, dotted, no head, from=5-4, to=4-6]
        \arrow[no head, from=4-1, to=3-1]
        \arrow[no head, from=4-3, to=3-3]
        \arrow[no head, from=4-6, to=3-6]
        \arrow[no head, from=4-8, to=3-8]
        \arrow[curve={height=12pt}, no head, from=1-4, to=3-1]
        \arrow[no head, from=1-4, to=3-3]
        \arrow[no head, from=1-4, to=3-6]
        \arrow[curve={height=-12pt}, no head, from=1-4, to=3-8]
      \end{tikzcd}\]
      \captionof{figure}{Path with hexagons}
      \label{fig:dfc_to_zpo_hasse_15}
    \end{minipage}
    \begin{minipage}{.5\textwidth}
      \[\begin{tikzcd}[column sep = 2.5pt]
        {d_0} && {d_1} && {d_{\varphi(i_0) - 1}} && {\hspace{-30pt} d_{\varphi(i_0)} = d_0} \\
        & {e_0} && {e_1} && {e_{\varphi(i_0)}}
        \arrow["\shortmid"{marking}, no head, from=2-6, to=1-7]
        \arrow[no head, from=2-6, to=1-5]
        \arrow[no head, from=2-4, to=1-3]
        \arrow[no head, from=2-2, to=1-1]
        \arrow["\shortmid"{marking}, no head, from=1-3, to=2-2]
        \arrow["\shortmid"{marking}, dotted, no head, from=2-4, to=1-5]
      \end{tikzcd}\]
      \captionof{figure}{Impossible path}
      \label{fig:dfc_to_zpo_hasse_16}
    \end{minipage}
  \end{figure}
  \noindent where for every $i\ge0$, $c_i$ is chosen as the unique element in $\Lambda_{n-k+1}$ such that $d_i\prec^{-} c_i$.
  
  Using the \hyperref[prop:DFC_hexagon]{hexagon property} in the rightmost hexagon yields a path in the tree structure of $\delta\big(\gamma^{k-2}\omega\big)$ (the same as in the proof of \hyperref[prop:DFC_pencil_linearity]{pencil linearity}), which must be from $c_p$ to $c_0$ because $d_p$ is not a target.
  On the other hand, there is a simple zig-zag from $c_p$ to $c_0$ obtained by the previous construction.
  But because of the uniqueness of simple zig-zag between two nodes of a tree, the path from $c_p$ to $c_0$ and the simple zig-zag from $c_0$ to $c_p$ constituted of triplets
  $$c_{\varphi(i),\,\sigma(i)}\succ^{\alpha_i}d_{\varphi(i),\,\sigma(i)}\prec^{-\alpha_i}c_{\varphi(i),\,\sigma(i)+1}$$
  for $0\le i < l$ must be the symmetric of each other.
  Hence it only contains such triplets with $\alpha_i = -$, and $e_i \prec^{+}d_{\varphi(i),\,\sigma(i)}$.
  In particular the last triplet $\cdot\succ^{+}d_0\prec^{-}c_0$ of the path from $c_p$ to $c_0$ corresponds to one of the triplets above.
  Suppose it is the $i_0$'th, then we have $e_{\varphi(i_0)} \prec^{+} d_{\varphi(i_0),\,\sigma(i_0)} = d_0$.
  Hence we have a path as in \ref{fig:dfc_to_zpo_hasse_16} above; this is absurd by \hyperref[defin:PFS]{strictness}, and ends the proof of the lemma.
\end{proof}

\begin{lem}\label{lem:DFC_backward_cell_path}
  If $e \prec^{-} d \prec^{-} c$, then there is a (unique) upper path from a source $e'$ of $\gamma(c)$ to $e$ as follows:
  
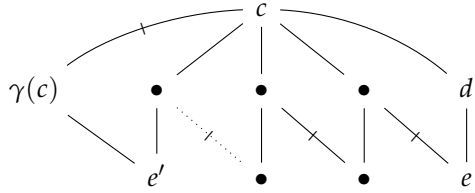
\[\begin{tikzcd}[row sep = 15pt]
    && c \\
    {\gamma(c)} & \bullet & \bullet & \bullet & d \\
    & {e'} & \bullet & \bullet & e
    \arrow["\shortmid"{marking}, no head, from=2-4, to=3-5]
    \arrow[no head, from=2-4, to=3-4]
    \arrow["\shortmid"{marking}, no head, from=2-3, to=3-4]
    \arrow[no head, from=2-3, to=3-3]
    \arrow["\shortmid"{marking}, dotted, no head, from=2-2, to=3-3]
    \arrow[no head, from=2-2, to=3-2]
    \arrow[no head, from=2-1, to=3-2]
    \arrow[no head, from=2-5, to=3-5]
    \arrow[no head, curve={height=-10pt}, from=1-3, to=2-5]
    \arrow[no head, from=1-3, to=2-4]
    \arrow[no head, from=1-3, to=2-3]
    \arrow[no head, from=1-3, to=2-2]
    \arrow["\shortmid"{marking}, no head, curve={height=10pt}, from=1-3, to=2-1]
  \end{tikzcd}\]
\end{lem}
\begin{proof}
  We keep completing half lozenges from right to left as above, until coming across the leftmost pattern.
\end{proof}

\begin{prop}\label{prop:DFC_disjointness_2}
  If $k>0$, and $d,\,d'\in C_k$ then it cannot be the case that $d<^{+}d'$ and $d<^{-}d'$.
\end{prop}
\begin{proof}
  Suppose that we have $d <^{+} d'$ and $d <^{-} d'$.
  Then by concatenating paths obtained as in \reflem{lem:DFC_backward_cell_path}, $d<^{+}d'$ yields a non-trivial upper path from $e$ to $\gamma(d)$ where $e\prec^{-}d'$ is a source of $d'$.
  On the other hand, $d<^{-}d'$, hence $\gamma(d)\le^{+}e'$ for some $e'$ a source of $d'$.
  Now by concatenating the path from $e$ to $\gamma(d)$ and the path from $\gamma(d)$ to $e'$, we obtain a non-trivial upper path from $e\prec^{-}\gamma(d')$ to $e'\prec^{-}\gamma(d')$, which is impossible by \reflem{lem:DFC_no_upper_path_in_cell}.
\end{proof}

\begin{prop}[disjointness]
  $C$ satisfies the axiom of \hyperref[defin:PFS]{disjointness}.
  That is: two elements $d$ and $d'$ cannot be comparable for both relations $<^{+}$ and $<^{-}$.
\end{prop}
\begin{proof}
  This is a consequence of \refprop{prop:DFC_disjointness_1} and \refprop{prop:DFC_disjointness_2}.
\end{proof}

\begin{thm}\label{thm:DFC_to_ZPO}
  The dendritic face complex $C$ is a \hyperref[defin:PFS]{positive opetope}.
\end{thm}
\begin{proof}
  This is a consequence of Propositions \ref{prop:DFC_glob}, \ref{prop:DFC_strictness_1}, \ref{prop:DFC_strictness_2}, \ref{prop:DFC_disjointness_1}, \ref{prop:DFC_disjointness_2}, \ref{prop:DFC_pencil_linearity} and \ref{prop:DFC_principal}.
\end{proof}

  
\section{From \textsc{Zawadowski}'s positive opetopes to dendritic face complexes}
\label{sec:zpo_to_dfc}

\para In this section, we consider an opetopic cardinal $S = \left((S_k)_{k\in\N},\,\gamma,\,\delta\right)$.
We aim to prove that $S$ gives rise to a dendritic face complex. (At a certain point, we will require $S$ to be principal in order to conclude.)
Recall from \refthm{thm:POP_pHg_equiv} that $S$ may be given the structure of a $POP$, we will then focus on proving the axioms of DFCs, in the following order: \hyperref[defin:DFC]{oriented thinness}, \hyperref[defin:DFC]{acyclicity} and the existence of a \hyperref[defin:DFC]{greatest element}.

\begin{prop}[Oriented Thinness]\label{prop:PFS_oriented_thinness}
  Let $S$ be an opetopic cardinal,
  $S$ seen as a positive-to-one poset satisfies the property of \hyperref[defin:DFC]{oriented thinness}.
\end{prop}
\begin{proof}
  We consider a chain $c \prec b \prec a$ in $S$.
  And we distinguish on the signs appearing in this relation.
  \begin{itemize}
    \item $\boxed{c \prec^{\beta} b \prec^{+} c}$ ~\\
      Because of \hyperref[defin:PFS]{globularity}, $\gamma\gamma(a) = \gamma\delta(a) - \delta\delta(a)$ and $\delta\gamma(a) = \delta\delta(a) - \gamma\delta(a)$.
      Hence there is some $b'$ as in \ref{fig:OT_lozenge_1} below.
      Suppose that there is another $b''$ such that $c\prec b'' \prec a$. \\
      Then still because of \hyperref[defin:PFS]{globularity}
      , either $c\prec^{-\beta}b''\prec^{+}a$ or $c\prec^{\beta}b''\prec^{-}a$. \\
      In the first case, $b'' = b$ by target uniqueness, hence $c\prec^{-} b$ and $c\prec^{+} b$ which is impossible.
      Hence only the second case may occur. \\
      In the second case, the point (2.) of the proposition $5.1$ in \cite{zawadowski2023positive} shows that necessarily $b'' = b'$, whence the uniqueness.
    
    \item $\boxed{c \prec^{+} b \prec^{-} c}$
      Lemma $4.1$ in \cite{zawadowski2023positive} gives $\gamma\delta(a) = \gamma\gamma(a) \sqcup \iota(a)$.
      Hence either $c\in\gamma\gamma(a)$ or $c\in\delta\delta(a)$ (those two cases are exclusive), and there is some $b'$ as in \ref{fig:OT_lozenge_2} below.
      Suppose that there is another $b''$ such that $c\prec b'' \prec a$. \\
      Then because of \hyperref[defin:PFS]{globularity} $c\notin\delta\gamma(a)$, hence there is some $\alpha'' \in\{+,\,-\}$ such that $c\prec^{\alpha''}b''\prec^{\alpha''}a$.
      And because the union $\gamma\delta(a) = \gamma\gamma(a) \sqcup \iota(a)$ is disjoint, $\alpha'' = \alpha'$.
      If $\alpha' = -$, we may conclude again by the point (2.) of the proposition $5.1$ in \cite{zawadowski2023positive} that $b' = b''$.
      And if $\alpha' = +$, then by uniqueness of the target, $b' = \gamma(a) = b''$.

    \item $\boxed{c \prec^{-} b \prec^{-} c}$
      Lemma $4.1$ in \cite{zawadowski2023positive} gives $\delta\delta(a) = \delta\gamma(a) \sqcup \iota(a)$.
      Hence either $c\in\delta\gamma(a)$ or $c\in\gamma\delta(a)$ (those two cases are exclusive), and there is some $b'$ as in \ref{fig:OT_lozenge_3} below.
      Suppose that there is another $b''$ such that $c\prec b'' \prec a$. \\
      Then because of \hyperref[defin:PFS]{globularity} $c\notin\gamma\gamma(a)$, hence there is some $\alpha'' \in\{+,\,-\}$ such that $c\prec^{-\alpha''}b''\prec^{\alpha''}a$.
      And because the union $\delta\delta(a) = \delta\gamma(a) \sqcup \iota(a)$ is disjoint, $\alpha'' = \alpha'$.
      If $\alpha' = -$, we may conclude again by the point (2.) of the proposition $5.1$ in \cite{zawadowski2023positive} that $b' = b''$.
      And if $\alpha' = +$, then by uniqueness of the target, $b' = \gamma(a) = b''$.
      \qedhere
  \end{itemize}
\end{proof}
\begin{figure}[H]
  \begin{minipage}{.33\textwidth}
    \[\begin{tikzcd}
      & a \\
      b && {b'} \\
      & c
      \arrow["\shortmid"{marking}, no head, from=1-2, to=2-1]
      \arrow["\beta"', no head, from=2-1, to=3-2]
      \arrow[no head, from=1-2, to=2-3]
      \arrow["{\beta}", no head, from=2-3, to=3-2]
    \end{tikzcd}\]
    \captionof{figure}{Lozenge 1}
    \label{fig:OT_lozenge_1}
  \end{minipage}%
  \begin{minipage}{.33\textwidth}
    \[\begin{tikzcd}
      & a \\
      b && {b'} \\
      & c
      \arrow[no head, from=1-2, to=2-1]
      \arrow["\shortmid"{marking}, no head, from=2-1, to=3-2]
      \arrow["{\alpha'}", no head, from=1-2, to=2-3]
      \arrow["{\alpha'}", no head, from=2-3, to=3-2]
    \end{tikzcd}\]
    \captionof{figure}{Lozenge 2}
    \label{fig:OT_lozenge_2}
  \end{minipage}%
  \begin{minipage}{.33\textwidth}
    \[\begin{tikzcd}
      & a \\
      b && {b'} \\
      & c
      \arrow[no head, from=1-2, to=2-1]
      \arrow[no head, from=2-1, to=3-2]
      \arrow["{\alpha'}", no head, from=1-2, to=2-3]
      \arrow["{-\alpha'}", no head, from=2-3, to=3-2]
    \end{tikzcd}\]
    \captionof{figure}{Lozenge 3}
    \label{fig:OT_lozenge_3}
  \end{minipage}
\end{figure}

\begin{prop}[Acyclicity]\label{prop:PFS_tree_structure}
  Let $S$ be an opetopic cardinal,
  $S$ satisfies the axiom of \hyperref[defin:DFC]{acyclicity}.
\end{prop}
\begin{proof}
  If $x \in S_1$, then $\delta(x)$ is a singleton because $\delta_0$ is functional.
  If $x \in S_{\ge1}$, then $\delta(x) \neq \emptyset$ because for all $k$, $\delta_k$ is total.
  There is no cycle as in \ref{fig:dfc_def_cycle} because of \hyperref[defin:PFS]{strictness}.
\end{proof}

\begin{prop}[Greatest element]\label{prop:principal_PFS_max}
  If $S$ is supposed to be principal, then $\hat{S}$ admits a greatest element $\omega$.
\end{prop}
\begin{proof}
  Because of \hyperref[defin:PFS]{principality}, $\vert \{ S_n \setminus \delta(S_{n+1})\} \vert = 1$.
  But $S_{n+1}$ is empty, hence $S_n$ is a singleton.
  The fact that its only element is indeed a greatest element is given by the point (1.) of Lemma $7.1$ in \cite{zawadowski2023positive}.
\end{proof}

\begin{thm}\label{thm:ZPO_to_DFC}
  Let $S$ be a positive opetope,
  $S$ is a dendritic face complex.
\end{thm}
\begin{proof}
  This is a consequence of \refprop{prop:PFS_oriented_thinness}, \refprop{prop:PFS_tree_structure} and \refprop{prop:principal_PFS_max}.
\end{proof}

\para Note that \refthm{thm:ZPO_to_DFC} is the converse of \refthm{thm:DFC_to_ZPO}.
 Because of our previous results, we may also state a converse to \refprop{prop:principal_PFS_max}:

\begin{thm}
  If an opetopic cardinal $S$ admits a greatest element, then it is principal.
  That is, it is a positive opetope.
\end{thm}
\begin{proof}
  Because of \refprop{prop:PFS_oriented_thinness} and \refprop{prop:PFS_tree_structure}, we know that $S$ may be seen as a DFC.
  We may then use the \refprop{prop:DFC_principal} to conclude.
\end{proof}

  
\section*{Conclusion and related works}
\label{sec:ccl}

We have shown in this paper a way to associate a (\textsc{Zawadowski}'s) positive opetope to any dendritic face complex, and vice versa.
Noting that these constructions extend to morphisms, we obtain two functors $F$ and $G$ as depicted below.
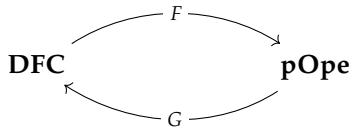
\[\begin{tikzcd}
	{\mathbf{DFC}} &&& {\mathbf{pOpe}}
	\arrow["F"{description}, curve={height=-24pt}, from=1-1, to=1-4]
	\arrow["G"{description}, curve={height=-24pt}, from=1-4, to=1-1]
\end{tikzcd}\]
where $\mathbf{DFC}$ and $\mathbf{pOpe}$ denote respectively the categories of dendritic face complexes and positive opetopes.
Since $G \circ F$ and $F \circ G$ leave the structure unchanged as proved in \refthm{thm:POP_pHg_equiv}, they form an equivalence of categories.

This result should be extended in a future paper, dealing with the equivalence with epiphytes (which will be defined then) and zoom complexes (see \cite{Kock2010}).

There is also a definition of opetopes, namely \textit{dendrotopes}, due to Thorsten \textsc{Palm} (see \cite{PalmArticle} for an introduction and \cite{PalmThesis} for a more complete description), which in many aspects is close to that of DFCs.
The author is convinced that there should be a functor from $\mathbf{DFC}$ to the category of \textsc{Palm}'s dendrotopes, although the details remain to be checked.

\section*{Acknowledgements}
\label{sec:ack}

Thanks to Amar \textsc{Hadzihasanović} for discussions about this article, and for ideas which helped develop the formalism of DFCs.
I also want to thank Marek \textsc{Zawadowski} for his inspiring and wide-ranging work on opetopes and their formalization as posets.
Finally, I would like to express a special thank you to the supervisor of my master's thesis, Pierre-Louis \textsc{Curien}, without whom this article would never have been.

  \printbibliography

@article{zawadowski2023positive,
  title={On positive opetopes, positive opetopic cardinals and positive opetopic set}, 
  author={Marek \textsc{Zawadowski}},
  year={2023},
  note = {\Arxiv{0708.2658}},
  primaryClass={math.GT}
}

@article{hadzihasanovic2019combinatorialtopological,
  title={A combinatorial-topological shape category for polygraphs}, 
  author={Amar \textsc{Hadzihasanovic}},
  year={2019},
  note = {\Arxiv{1806.103538}},
  primaryClass={math.CT}
}

@article{Kock2010,
	doi = {10.1016/j.aim.2010.02.012},
	year = 2010,
	month = {aug},
	publisher = {Elsevier {BV}},
	volume = {224},
	number = {6},
	pages = {2690--2737},
	author = {Joachim \textsc{Kock} and Andr{\'{e}} \textsc{Joyal} and Michael \textsc{Batanin} and Jean-Fran{\c{c}}ois \textsc{Mascari}},
	title = {Polynomial functors and opetopes},
	journal = {Advances in Mathematics},
  note = {\Arxiv{0706.1033}}
}

@article{Curien2022,
  note = {\url{https://higher-structures.math.cas.cz/api/files/issues/Vol6Iss1/CurHoTMim}},
  year = 2022,
  title = {Type theoretical approaches to opetopes},
  volume = {6},
  pages = {80--181},
  author = {C{\'{e}}dric \textsc{Ho Thanh} and Pierre-Louis \textsc{Curien} and Samuel \textsc{Mimram}},
  journal = {Higher Structures}
}

@article{baez1997higherdimensional,
  title={Higher-Dimensional Algebra III: n-Categories and the Algebra of Opetopes}, 
  author={John C. \textsc{Baez} and James \textsc{Dolan}},
  journal={Advances in Mathematics},
  doi={135(2):145–206},
  year={1998},
  note = {\Arxiv{q-alg/9702014}},
  archivePrefix={arXiv},
  primaryClass={q-alg}
}

@book{leinster2003higher,
      title={Higher Operads, Higher Categories}, 
      author={Tom \textsc{Leinster}},
      year={2003},
      primaryClass={math.CT},
      note = {\Arxiv{math/0305049}}
}

@misc{Ara2010,
  title = "Sur les infini-groupoïdes de Grothendieck et une variante infini-catégorique",
  author = "Dimitri \textsc{Ara}",
  year = "2010",
  pages = "167",
  note = "Thèse de doctorat dirigée par Maltsiniotis, Georges Mathématiques Paris 7 2010",
}

@article{PalmArticle,
  author          = {Thorsten \textsc{Palm}},
  title           = {Dendrotopic Sets},
  journal         = {Fields Institute Communications},
  title           = {Galois Theory, Hopf Algebras, and Semiabelian Categories},
  volume          = {43},
  year            = {2004},
  doi             = {https://doi.org/10.1090/fic/043}
}

@misc{PalmThesis,
  title = "Dendrotopic sets for weak infinity-categories",
  author = "Thorsten \textsc{Palm}",
  year = "2003",
  note = "Ph.D. - York University",
}

\end{document}